\documentclass[a4paper,11pt]{article}
\usepackage[top=18mm,left=22mm]{geometry}\textwidth17cm\textheight24.2cm
\usepackage{amsmath,amssymb,amsthm}
\usepackage{mathtools}
\usepackage[position=top]{subfig}
\usepackage{graphicx}
\usepackage[utf8]{inputenc}
\usepackage[]{graphicx} 
\usepackage{pstricks}
\usepackage{paralist}
\usepackage{float}
\usepackage[colorlinks=true]{hyperref}

\def\bmip{\begin{minipage}{\textwidth}}\def\emip{\end{minipage}}
\def\huga#1{\begin{gather} #1 \end{gather}}

\def\hual#1{\begin{align} #1 \end{align}}

\newcommand{\R}{{\mathbb R}}
\newcommand{\C}{{\mathbb C}}

\newcommand{\Ha}{{\cal H}}
\newcommand{\argmax}{\operatornamewithlimits{argmax}}

\def\ig{\includegraphics}

\def\CH{{\cal H}}

\def\del{\delta}
\def\ga{\gamma}\def\uh{\hat{u}}
\def\ds{\displaystyle}\def\al{\alpha}
\def\pa{{\partial}}
\def\Om{\Omega}
\def\dd{\, {\rm d}}\def\er{{\rm e}}
\newcommand{\bi}{\begin{itemize}}\newcommand{\ei}{\end{itemize}}
\newcommand{\ben}{\begin{enumerate}}\newcommand{\een}{\end{enumerate}}
\newcommand{\bci}{\begin{compactitem}}\newcommand{\eci}{\end{compactitem}}
\newcommand{\bcen}{\begin{compactenum}}\newcommand{\ecen}{\end{compactenum}}
\newcommand{\bce}{\begin{center}}\newcommand{\ece}{\end{center}}
\newcommand{\reff}[1]{(\ref{#1})}
\newcommand{\ov}[1]{{\overline {#1}}}
\newcommand{\spr}[1]{\left\langle #1 \right\rangle}

\def\aqui{\Leftrightarrow}
\def\ra{\rightarrow}

\newcommand{\barr}{\begin{array}}\newcommand{\earr}{\end{array}}
\newcommand{\bpm}{\begin{pmatrix}}\newcommand{\epm}{\end{pmatrix}}
\newcommand{\bsm}{\left(\begin{smallmatrix}}
\newcommand{\esm}{\end{smallmatrix}\right)}
\newcommand{\ba}{\begin{array}}\newcommand{\ea}{\end{array}}
\def\PDE{{PDE}}\def\PDEs{{\PDE s}}
\def\ODE{{ODE}}\def\ODEs{{\ODE s}}
\def\OCMAT{{\texttt{OCMat}}}

\def\OCPRO{{optimal control problem}} 
\def\OCPROs{{\OCPRO s}}
\def\PDETOPATH{{\texttt{pde2path}}}

\def\defin{{:=}} 
\def\cref#1{(\ref{#1})}
\def\crefsys#1{(\ref{#1})}\def\labelcref#1{(\ref{#1})}
\def\D{\operatorname{d}}

\def\E{\operatorname{e}}
\newcommand{\Dt}[1][t]{\D\!#1} %
\newcommand{\Dx}[1][x]{\D\!#1} 
\newcommand{\asymmat}{\Psi}

\def\PMAXP{Pontryagin's Maximum Principle}

\newcommand{\httpTOM}{\url{http://www.dm.uniba.it/~mazzia/bvp/index.html}}

\newcommand\css{{CSS}}
\newcommand\oss{{OSS}}

\newcommand\hcssnspp{{PCSS$^-$}}

\newcommand\eigval{\xi}

\newcommand\defect{\operatorname{d}}

\renewcommand\Re{\operatorname{Re}}
\newcommand*\evalat[1]{\left.#1\right\rvert}

\def\pdep{{\tt pde2path}}
\def\ph{\hat p}\def\qh{\hat q}
\def\Jloc{J_{\rm local}}\def\ddt{\frac{{\rm d} }{{\rm d}t}}
\newtheorem{theorem}{Theorem}[section]

\newtheorem{cor}[theorem]{Corollary}
\newtheorem{definition}[theorem]{Definition}
\newtheorem{remark}[theorem]{Remark}

\def\brem{\begin{remark}}\def\erem{\end{remark}}

\def\red{}
\def\hoss{POSS}\def\foss{FOSS}\def\fcss{FCSS}\def\hcss{PCSS}
\def\phoss{\hat P_{{\rm PS}}}
\def\phossnspp{\hat P_{{\rm PS}^-}}
\def\pfossone{\hat P_{{\rm FSC}}}

\def\pfosstwo{\hat P_{{\rm FSM}}}

\def\uhat{\hat{u}}\def\uh{\hat{u}}
\def\mlab{{\tt Matlab}} 
\def\eex{\hfill\mbox{$\rfloor$}}
\allowdisplaybreaks
\def\citet{\cite}\def\citep{\cite}
\begin{document}
\title{Optimal management and spatial patterns in a distributed shallow lake 
model}
\author{D. Gra\ss$^1$, H. Uecker$^2$\\ \small
$^1$ ORCOS, Institute of Mathematical Methods in Economics, Vienna University of Technology,\\ 
\small
A-1040 Vienna, Austria, dieter.grass@tuwien.ac.at\\
\small
$^2$ Institut f\"ur Mathematik, Universit\"at Oldenburg, D26111 Oldenburg, hannes.uecker@uni-oldenburg.de}

\normalsize

\maketitle

\begin{abstract} We present a numerical framework to treat 
infinite time horizon spatially distributed optimal control problems 
via the associated canonical system derived by  \PMAXP. 
The basic idea is to consider the  
canonical system in two steps. First we perform a bifurcation 
\mbox{analysis} of canonical steady states using the 
continuation and bifurcation package \pdep, yielding 
a number of so called flat and patterned canonical steady states. 
In a second step we link \pdep\ to the two point boundary value problem 
solver {\tt TOM} to study time dependent canonical system 
paths to steady states having the so called saddle point property. 
As an example we consider a shallow lake model with diffusion. 
\end{abstract}


\section{Introduction}
\label{sec:introduction}
In \cite{BX08} the authors consider economically motivated deterministic
optimal control (OC) problems with an infinite time horizon and a continuum 
of spatial sites
over which the state variable can diffuse. Using 
\PMAXP\ (see \S\ref{slsec2} for general background and references on 
OC problems in a PDE setting) they derive 
the associated canonical system and show the remarkable result 
that under certain conditions on the Hamiltonian 
there occurs a Turing like bifurcation 
from flat to patterned steady states of the canonical system, 
and call this phenomenon optimal diffusion-induced instability (ODI). 
Here we present a numerical 
framework to (a) study such ODI bifurcations of canonical steady states 
(CSS) numerically in a simple way, and (b) study their optimality by calculating 
and evaluating time--dependent paths to and from such CSS. 
As an example we use one of the three examples presented in 
\cite{BX08}, namely a version of
the, in the field of ecological economics, well-known shallow lake
optimal control (SLOC) model, 
cf.~\cite{scheffer1998,maeleretal2003,carpenterbrock2004}. 

We use the acronyms FCSS and PCSS 
for flat and patterned canonical steady states, and 
similarly FOSS and POSS for {\em optimal} canonical steady states, 
and summarize FOSS and POSS as OSS. 
The SLOC model has up to three (branches of) FCSS in relevant parameter regimes, 
and in these regimes we also find a large number of (branches of) PCSS.  
In this 
situation of multiple CSS a local stability 
analysis at a given CSS is in general not sufficient. 
Here, local stability analysis means that the
stationary canonical system is analyzed, analogous to a steady-state
analysis of the canonical system derived from an \OCPRO{} without
spatial diffusion. It is well known and shown for many models that the
appearance of multiple CSS (even if
these steady-states are saddle-points) does not necessarily imply the
appearance of multiple steady-states in the optimal system, 
cf.~\cite{grassetal2008,kiselevawagener2008a}. 
The reason is that there can exist a non--constant 
extremal solution, i.e.~a canonical path connecting the state values of a 
given CSS to some other CSS, and yielding 
a higher objective value. Therefore, to study whether a \css{} 
corresponds to an \oss, the values of the associated stable paths 
also have to be considered. 

Here we numerically compute the bifurcation behavior of 
\fcss{} and \hcss{} for the SLOC model in some detail, and 
study their optimality by evaluation of their objective 
values $J$ and comparison to time-dependent canonical paths. 
Such a global analysis is 
inevitable and has to accompany the local stability analysis.
Since in general the pertinent \ODEs{} or \PDEs\ cannot be solved
analytically we have to use numerical methods for the 
calculation of \fcss{} and in particular \hcss, and for the calculation 
of $t$--dependent canonical paths. 
For the steady state problem we use the 
continuation and bifurcation software \PDETOPATH\ \cite{p2pure}, based on a 
spatial finite element method (FEM) discretization, which we then 
combine with the boundary value problem (BVP) solver TOM to obtain 
canonical paths. 

A standard reference on ecological economics or 
``Bioeconomics'' is \cite{Cl90}, which also 
contains a very readable 
account, and applications, of \PMAXP{} in the context of 
ODE models, while \cite{grassetal2008} focuses more on 
socio--economical ODE model applications.  
Besides in \cite{BX08}, and in \cite{BX10}, PDE models roughly similar to 
our diffusive SLOC model are considered in, e.g., 
\cite{LW07,AAC11, DHM12, ACKT13, Apre14}, partly including 
numerical simulations.  
However, these works are in a finite time horizon setting, 
and with control constraints, which 
altogether gives a rather different setting from the one considered here. 
See Remark \ref{lrem} for further comments. 

In \S\ref{msec} we present the SLOC model. 
To give some background, in \S\ref{slsec1} we briefly present the 
0D (ODE) version, and 
some basic concepts of optimal control, in particular \PMAXP. 
In \S\ref{slsec2} we turn to the distributed case, explain the 
associated canonical system, and relate our method 
to other approaches to PDE OC problems. 
In \S\ref{nsec} we explain the numerics 
to first compute the bifurcation diagram of canonical steady states, and 
then to solve the BVP in $t$. 
We mostly focus on one spatial dimension (1D), but also give a short 
outlook on the 2D case. {\red It turns out that in the  parameter 
regimes studied here the PCSS are not optimal, but nevertheless 
they play a relevant role. Moreover, calculating {\em optimal} canonical 
paths to FCSS yields interesting and to some extent counter--intuitive 
information about the optimal control of the distributed SLOC model. 

In \S\ref{dsec} we close with a short summary and discussion, and 
Appendix \ref{sec:spppdes} contains remarks about the saddle point 
property for CSS in a PDE setting, starting on the discretized level.

Our software, including demo files 
and a manual to run some of the simulations in this paper, 
can be downloaded from 
\url{www.staff.uni-oldenburg.de/hannes.uecker/pde2path}. 
In fact, the present paper is the first in a series of 
four related works, the other three being \cite{U15p2, grass2014, U15}: 
\cite{U15p2} contains a Quickstart guide and implementation details 
of the add on package {\tt p2poc} to {\tt pde2path} used for the 
computations in this paper. 
Thus, the reader interested in these details should read 
(parts of) \cite{U15p2} in parallel.  Next, \cite{grass2014} 
explains the usage of {\tt OCMat} \url{orcos.tuwien.ac.at/research/ocmat_software/} to study 1D distributed OC problems 
based on spatial finite difference approximations, with the 
same SLOC model as in the present paper as an example, and thus obtaining 
comparable results, but also studying a second parameter regime.  
Finally, in \cite{U15} we 
apply {\tt p2poc} to an OC problem for a 
reaction--diffusion {\em system} modeling a vegetation-water-grazing 
interaction. In contrast to the SLOC model studied here,  this yields 
dominant patterned optimal steady states in wide parameter 
regimes, and thus interesting new results on spatial patterns in 
optimal harvesting. 
}

\section{The model, and background from optimal control}\label{msec} %
\subsection{The shallow lake model without diffusion}
\label{slsec1}
A well known non--distributed or 0D version of the SLOC model,
see e.g.~\cite{wagener2003}, can be formulated in dimensionless form as
\begin{subequations}
  \label{sl_model}
  \begin{align}
    V(P_0)&:=\max_{k(\cdot)}J(P_0,k(\cdot)), \qquad 
J(P_0,k(\cdot)):=\int_0^\infty\E^{-r t} J_c(P(t),k(t))
\!\Dt, \label{sl_model_obj2}\\
\text{ where}\quad J_c(P,k)&=\ln k-\ga P^2\\
\intertext{is the current value objective function, and $P$ fulfills 
the ODE initial value problem} 
 \dot P(t)&=k(t)-bP(t)+\frac{P(t)^2}{1+P(t)^2}\label{sl_model_dyn},\qquad 
    P(0)=P_0\ge0.
  \end{align}
\end{subequations} 
Here $r,\ga,b>0$ are parameters, $P=P(t)$ is the phosphorus 
contamination of the lake, which we want to keep low for ecological
reasons, and $k=k(t)$ is the phosphate load, for instance from 
fertilizers used by farmers, which farmers want high for economic
reasons. The objective function consists of the concave
increasing function $\ln k$, and the concave decreasing function
$-\ga P^2$;  $b$ is the phosphorus degradation 
rate in the lake, and $r$ is the discount rate. The discounted 
time integral in \reff{sl_model_obj2} is typical for economic (or socio-political) 
problems, where ``profits now'' weight more than mid or far 
future profits. More specifically, $r$ often 
corresponds to a long-term investment rate. 
We focus on the parameter choice 
\hual{ r=0.03,\quad \ga=0.5,\quad b\in (0.5,0.8) \text{ (primary
    bifurcation parameter),}\label{parsel}
}
and for the distributed case we shall additionally fix the  
diffusion parameter to $D=0.5$. 

 The $\max$ in (\ref{sl_model_obj2}) runs over all {\em admissible} 
controls $k$ and (associated) states $P$; for $k$ we can take the space 
$C^0_b([0,\infty),\R)$, and for $P$ the space $C^1_b([0,\infty),\R)$. 
In fact, we naturally have $k(t)> 0$ for all $t$ as $J_c(k,P)\ra -\infty$ 
as $k\searrow 0$, and then \reff{sl_model_dyn} implies that $P(t)>0$ for all 
$t>0$. On the other hand, see Remark \ref{lrem} for comments on the case  
of state or control constraints. 

By \PMAXP{} \cite{pontryaginetal1962}, see also, e.g., 
\cite{grassetal2008}, 
an optimal solution $(P^*(\cdot),k^*(\cdot))$ has to satisfy
the first order optimality conditions 
\begin{subequations}
  \label{eq:hammax}
  \begin{align}
    &k^*(t)=\argmax_{k}\Ha(P(t),k,q(t),q_0)\quad\text{for
      almost all}\quad t\ge0,\label{eq:hammax1} \shortintertext{with
      the local current value Hamiltonian function}
    &\Ha(P,k,q,q_0):=q_0J_c(P,k)
+q\left(k-bP+\frac{P^2}{1+P^2}\right).\label{eq:hammax2}
  \end{align}
\end{subequations}
The state $P(\cdot)$ and costate $q(\cdot)$ paths are solutions of
the canonical system\footnote{It can be proved that the problem is
  normal, i.e. $q_0>0$, and hence w.l.o.g. $q_0=1$ can be
  assumed and is therefore subsequently omitted.}
\begin{subequations}
  \label{cansys1}
  \begin{align}
    &\dot P(t)=\pa_q\Ha(P(t),k^*(t),q(t))=k^*(t)-bP(t)+\frac{P(t)^2}{1+P(t)^2},
    \label{sl_model_cansys1}\\
    &\dot q(t)=rq(t)-\pa_P\Ha(P(t),k^*(t),q(t))=2\ga P(t)+q(t)
    \left(r+b-\frac{2P(t)}{\left(1+P(t)^2\right)^2}\right),
    \label{sl_model_cansys2}
    \shortintertext{with} &P(0)=P_0>0,\notag
  \end{align}
\end{subequations}
additionally satisfying the transversality condition
\begin{equation}
  \label{tcond}
  \lim_{t\to\infty}\er^{-rt}q(t)=0\quad\text{if}\quad\liminf_{t\to\infty}P^*(t)>0.
\end{equation}
A solution $(P(\cdot),q(\cdot))$ of the canonical system 
\reff{cansys1} is called a 
\emph{canonical path}, 
and a steady state of \reff{cansys1} is called a 
\emph{canonical steady state (\css{})}.
Due to the strict concavity and continuous differentiability of the
Hamiltonian function with respect to the control $k$, and the absence of 
control constraints, the solution of
\reff{eq:hammax1} is given by
\begin{equation}
  \label{eq:solham}
  \pa_k\CH(P(t),k(t),q(t))=0\quad\text{which yields}\quad k^*(t)=-\frac{1}{q(t)}.
\end{equation}
Consequently, for $(P(\cdot),q(\cdot))$ a canonical path, i.e., 
a solution of the canonical system, 
with a slight abuse of notation we also call $(P,k)$ with $k=-1/q$ a 
canonical path. In particular, if $(\hat P, \hat q)$ is a CSS, so is 
$(\hat P,\hat k)$. 
Canonical paths yield candidates for optimal solutions, 
defined as follows: 
\begin{definition}
  \label{def:optsys}
  $(P^*(\cdot),k^*(\cdot,P_0))$ is called an optimal solution of 
\reff{sl_model} if for every admissible $k(\cdot)$ and 
associated $P(\cdot)$ we have
	\begin{equation*}
			J(P_0,k(\cdot))\le J(P_0,k^*(\cdot,P_0))=V(P_0).
	\end{equation*}
	Then $k^*(\cdot,P_0)$ is called an \emph{optimal control}, $P^*(\cdot)$ the corresponding \emph{optimal (state) path}, and 
  \begin{equation}
    \label{eq:optsys}
    \dot P(t)=k^*(t,P_0)-bP(t)+\frac{P(t)^2}{1+P(t)^2}
  \end{equation}
  is called the \emph{optimal ODE}. A constant solution $(P^*(\cdot),k^*(\cdot,P_0))\equiv(\hat P,\hat k(\hat P))$ of \reff{eq:optsys} is called an \emph{optimal steady state (\oss)}.
\end{definition}
It turns out that the long-run behavior of an optimal 
solution $(P^*(\cdot),k^*(\cdot))$ can be characterized 
completely, see, e.g., \cite{wagener2003}. Each optimal solution 
converges to an OSS, and depending on the parameters 
\reff{cansys1} can have $I=1,2,3$ \css{} $(\hat P,\hat q)_i$, 
$i=1,\ldots,I$.\footnote{ 
cf., e.g., the FCSS branches in Fig.~\ref{f1} for our specific parameter choice.}
For simplicity omitting the non-generic case $I=2$, if $I=1$ then the 
unique CSS is a globally stable OSS, while for $I=3$ two CSS are locally 
stable OSS, and the third is unstable. Here a OSS $(\hat P,\hat q)$ is 
called globally (locally) stable if for each $P(0)$ (in a neighborhood 
of $\hat P$) the associated optimal path converges to $(\hat P,\hat q)$; see 
\cite{kiselevawagener2008a,kiseleva2011} for a detailed discussion.

Setting $u:=(P,q)$ and letting $\hat{u}$ be a steady state of
\reff{cansys1}, the problem now is to compute a path, or 
``connecting orbit'', with $P(0)=P_0$ and 
$\ds\lim_{t\to\infty}u(t)=\hat u$. One standard approach, see, 
e.g.~\cite{lentinikeller1980, auto, BPS01} and in particular 
\cite[Chapter 7]{grassetal2008}, is to treat
\reff{cansys1} on a finite time interval $[0,T]$ and to require
$u(T)\in W_s(\hat u)$, where $W_s(\hat u)$ is the local stable
manifold of $\hat u$. In practice we approximate $W_s(\hat u)$ by the
stable eigenspace $E_s(\hat u)$, and thus require 
\huga{\label{wscond}
  u(T)\in E_s(\hat u)\text{ and close to }\hat u.  
}  
To obtain a well 
defined two point boundary value problem we then need $\dim E_s(\hat
u)=1$.

More generally, if the state variable is an $n$--dimensional vector and 
thus the canonical system is a system of $2n$ ODEs, for arbitrary $P_0$ 
we need that $\hat{u}$ has the saddle point property, defined 
as follows. 
\begin{definition}[Saddle Point Property]
\label{def:spp}
A \css{} $\uh\in\R^{2n}$ with 
\begin{equation} n_s:=\dim E_s(\uh)=n \label{spp} \end{equation}
is called a \css{} with the \emph{saddle point property (SPP)}. 
The number $\defect(\uh)\defin n_s-n$ 
is called the \emph{defect} of $\uh$, and a \css{} 
with defect $\uh<0$ is called \emph{defective}. 
\end{definition}
For ODE problems like \reff{cansys1}, 
given $\uh$ with the SPP, and some initial state $P_0$, 
 canonical paths 
connecting $P_0$ and $\uh$ can now be computed using, e.g., 
{\tt OCMat}. We now generalize this to distributed problems, and thus 
in \S\ref{nsec} explain further details on that level. 

\subsection{The shallow lake model with diffusion}
\label{slsec2}
Following \cite{BX08} we consider the shallow lake model with diffusion 
in a domain $\Om\subset\R^d$, $d=1,2$, i.e.,
\begin{subequations} 
\label{sldiff1}
  \begin{align}
    &V(P_0(\cdot)):=\max_{k(\cdot,\cdot)}J(P_0(\cdot),k(\cdot,\cdot)), \qquad 
    J(P_0(\cdot),k(\cdot,\cdot)):=\int_0^\infty\E^{-r t}
J_{ca}(P(t),k(t))\Dt,\label{sldiffusion_model_obj2}\\
\text{where }& J_{ca}(P(\cdot,t),k(\cdot,t))=\frac 1 {|\Om|}
\int_\Om J_c(P(x,t),k(x,t))\!\Dx
\quad\text{($J_c(P,k)=\ln k-\ga P^2$ as before)}
\\
\intertext{is the spatially averaged current value objective 
function, and $P$ fulfills the initial boundary value problem 
}
 &\pa_t P(x,t)=k(x,t)-bP(x,t)+\frac{P(x,t)^2}{1+P(x,t)^2}+D\Delta P(x,t),\label{sldiffusion_model_dyn}\\
    &\evalat{\pa_\nu P(x,t)}_{\pa\Omega}=0,\label{sldiffusion_model_bc}\qquad 
    \evalat{P(x,t)}_{t=0}=P_0(x),\quad
    x\in\Om\subset\R^d,
  \end{align}
\end{subequations}
where $\Delta=\pa_{x_1}^2+\ldots+\pa_{x_d}^2$, and $\nu$ is the outer 
normal of $\Om$. We normalize $J_{ca}$ by the domain size $|\Om|$ 
for easy comparison between the 0D, 1D, and 2D cases, and, more 
generally, between different domains. We mostly focus on 
$\Om=(-L,L)$ a real interval, but also give an outlook to 2D. In 2D 
the model is somewhat less intuitive,  as 
a controlled phosphate dumping in the ``middle'' of the lake 
from farming appears difficult to motivate, 
and thus in 2D we rather think of \reff{sldiff1} as a general pollution model. 
Instead of the periodic BC in 1D in \cite{BX08}
we require Neumann (zero flux)
boundary conditions (BC), which from a
modeling point of view we find more natural.  

Introducing the costate $q:\Om\times(0,\infty)\ra \R^{N}$ and the 
(local current value) Hamiltonian
\hual{
\CH&=\CH(P,q,k)=J_c(v,k)+q \bigl[k-bP+\frac{P^2}{1+P^2}+D\Delta P\bigr],
}
by Pontryagin's Maximum Principle for $\tilde{\CH}=
\int_0^\infty \er^{-r t} \ov{\CH}(t)\dd t$ with the spatial integral 
\huga{\label{fullH} 
\ov{\CH}(t)=\int_\Om \CH(P(x,t),p(x,t),k(x,t))
\dd x, 
}
the canonical system for \reff{sldiff1} becomes 
\begin{subequations}
  \label{slcan2}
  \begin{align}
&    {\pa_t} P(x,t)=[\pa_qH](x,t)=k(x,t)-bP(x,t)+\frac{P(x,t)^2}{1+P(x,t)^2}
+D\Delta P(x,t),\label{sldiffusion_model_cansys1}\\
 &   \pa_t q(x,t)=rq(x,t)-[\pa_P H](x,t)=2\ga P(x,t)+q(x,t)\left(r+b-\frac{2P(x,t)}{\left(1+P(x,t)^2\right)^2}\right)-D\Delta q(x,t),\label{sldiffusion_model_cansys2}\\
&    \evalat{\pa_\nu P(x,t)}_{\pa\Omega}=0,\quad
\evalat{\pa_\nu q(x,t)}_{\pa\Omega}=0,\label{sldiffusion_model_cansys3}\qquad 
    \evalat{P(x,t)}_{t=0}=P_0(x),\quad x\in\Omega, 
\intertext{where $k=\argmax_{\tilde{k}}\CH(P,q,\tilde{k})$, 
which similar to \reff{eq:solham} is obtained from}
&\pa_k \CH(P,q,k)=0\quad\aqui\quad k(x,t)=-\frac{1}{q(x,t)} 
 \end{align}
\end{subequations}
The costate $q$ also fulfills zero flux BC, and 
derivatives like $\pa_P \CH$ etc are taken 
variationally, i.e., for $\ov{\CH}$. For instance, for 
$\Phi(P,q):=q\Delta P$ we have $\ov{\Phi}(P,q)
=\int_\Om q\Delta P\dd x
=\int_\Om (\Delta q)P\dd x$ by Gau\ss' theorem, hence 
$\delta_P \ov{\Phi}(P,q)[h]=\int (\Delta q) h\dd x$, and 
by the Riesz representation theorem we 
identify $\delta_P \ov{\Phi}(P,q)$ and hence $\pa_P\Phi(P,q)$ 
with the multiplier $\Delta q$. 

Finally, we have the limiting intertemporal transversality condition 
(see Remark \ref{lrem} below) 
\begin{equation}
\label{tcond2}
	\lim_{t\ra\infty}\er^{-rt}\int_\Omega q(x,t)P(x,t)\!\Dx=0.
\end{equation}

Analogous to Def.~\ref{def:optsys} we define
\begin{definition}
  \label{def:diffoptsys}
  Let $(P^*(\cdot,\cdot),k^*(\cdot,\cdot,P_0))$ be an optimal 
solution of problem \reff{sldiff1}, i.e.~for every 
admissible $k(\cdot,\cdot)$ and associated $P(\cdot,\cdot)$ we have
	\begin{equation*} 
J(P_0,k(\cdot,\cdot))\le J(P_0,k^*(\cdot,\cdot,P_0))=V(P_0).
	\end{equation*}
	Then $k^*(\cdot,\cdot,P_0)$ is called a (distributed) 
\emph{optimal control}, $P^*(\cdot,\cdot)$ is called the associated 
distributed \emph{optimal (state) path}, and 
  \begin{align}
    \label{eq:optpdesys}
    &\pa_t P(x,t)=k^*(x,t,P_0(x))-bP(x,t)+\frac{P(x,t)^2}{1+P(x,t)^2}+D\Delta P(x,t), \quad  \evalat{\pa_\nu P(x,t)}_{\pa\Omega}=0
  \end{align}
  is called the \emph{optimal PDE}. Again with a slight abuse of notation, 
$(P^*,k^*)$ is also called an optimal path, and 
an optimal stationary solution $(\hat P(\cdot),\hat k(\cdot))$ of \reff{slcan2} 
is called an \emph{\oss} \emph{(optimal steady state)}. 
If $\hat P(\cdot)\equiv\hat P$ then the optimal steady state is called 
a \emph{\foss{} (flat optimal steady state)}, otherwise it is called a 
\emph{\hoss{} (patterned optimal steady state)}.
\end{definition}
For a CSS $\uh(\cdot)$ we additionally introduce the acronyms 
\emph{\fcss{}} for a flat canonical steady state, 
i.e $\uh(\cdot)\equiv\uh$, and \emph{\hcss{} 
(patterned canonical steady state)} otherwise. Obviously, 
the FCSS correspond precisely to the 0D CSS from \S\ref{slsec1}. 
It was already indicated 
in \cite{BX08} that \reff{slcan2} can additionally have 
PCSS arising from Turing like
bifurcations.  Thus, we first calculate bifurcation diagrams for
\reff{slcan2}, in 1D and 2D, recovering the up to 3 branches 
of \fcss\ from \S\ref{slsec1}, and many branches of PCSS. Next, 
analogous to the 0D case, we expect a solution of \reff{sldiff1} to 
converge to some CSS. Thus, we only 
consider solutions $u(\cdot,\cdot)$ 
of \reff{slcan2} with $\lim_{t\to\infty} u(\cdot,t)=\hat u(\cdot)$, 
where $\hat u(\cdot)$ is a CSS. For 
$\hat u$ we then also need a version of the SPP. 
However, $E_s(\uh)$ (and $W_s(\uh)$) and 
$E_u(\uh)$ (and $W_u(\uh)$) are infinite
dimensional. 
We circumvent this problem by requiring
\reff{wscond} and the saddle point property after a spatial 
discretization, which
turns \reff{slcan2} into a (very large) systems of ODEs again. 
See App.~\ref{sec:spppdes} for further discussion. 

{\red \brem\label{lrem}{\rm 
(a) A strict mathematical proof of \PMAXP{} for diffusion processes over 
an infinite time horizon is still missing, specifically for the 
transversality condition \reff{tcond2}, which also in \cite{BX08} is discussed 
only on a heuristic base. 
Thus, at the moment we apply \PMAXP{} in an ad hoc sense. 
We specifically assume, 
based on the results for the 0D shallow lake model, that canonical paths 
converge to CSS, and therefore make no use of the ``critical'' transversality 
condition \reff{tcond2}. In any case, a particular feature of the 
canonical system for diffusion processes is 
the anti--diffusion in the co--states, cf.~(\ref{slcan2}b), which 
makes the canonical system ill--posed as an initial value problem. 

(b) For background on OC in a PDE setting see also for instance 
\cite{Tr10} and the references therein, or specifically 
\cite{RZ99, RZ99b,DHM12,ACKT13} and 
\cite[Chapter5]{AAC11} for \PMAXP{} 
for OC problems for semilinear parabolic state evolutions. 
However, these works are in a 
finite time horizon setting, and often the objective function 
is linear in the control and there are control constraints, e.g., 
$k(x,t)\in K$ with some bounded interval $K$. Therefore $k$ is not obtained from 
the analogue of (\ref{slcan2}d), but rather takes the values from $\pa K$, which 
is usually called bang--bang control.  
In, e.g., \cite{Neu03} and \cite{DL09}, stationary 
spatial OC problems for a fishery model with (active) control constraints 
are considered, including 
numerical simulations, which correspond to our calculation of 
canonical steady states for our SLOC model. 
Here we do not (yet) consider explicit control or state constraints, 
and have an objective function strictly concave in the control, and thus we 
have a rather different setting than the above works. 

(c) We summarize that we do not aim at new theoretical results, but 
rather consider \reff{slcan2} after a spatial discretization as a 
(large) ODE problem, to which we apply the ``connecting orbit method''. 
Importantly, this means that we take a 
broader perspective than aiming at computing just one optimal control, 
given an initial condition $P_0$, which without further information 
is an ill-posed problem anyway. 
Instead, our method aims to give a somewhat global 
picture by identifying the pertinent CSS and their respective domains 
of attraction. \eex 
}
\erem }

\section{Numerical algorithm, and results}
\label{nsec}
The general idea is to use a method of lines discretization 
of \reff{slcan2}, i.e., to approximate 
\huga{\label{mola} u(x,t):=(P(x,t),q(x,t))=\sum_{i=1}^{2n}
  u_i(t)\phi_i(x), 
} 
where $(\phi_i)_{i=1,\ldots,2n}$ spans a subspace
$X_n$ of the phase space $X$ of \reff{slcan2}, e.g., here
$X=[H^1(\Om)]^2$, and $(u_i,u_{n+i})$, $i=1,\ldots,n$,  
are the expansion coefficients of $P,q$, respectively.  
This converts \reff{slcan2} into a (high dimensional) ODE
\begin{subequations}
  \huga{\label{zivp} \dot u(t)=- G(u(t)),\quad u_i(0)=u_{0,i},\ \ 
    i=1,\ldots,n, } for the coefficient vector
  $u=(u_i)_{i=1,\ldots,2n}$, where we have initial data for exactly
  half of the expansion coefficients.  We choose a truncation time $T$ 
and augment \reff{zivp} with the
  approximate transversality condition 
\huga{\label{asmf-cond} u(T)\in E_s(\hat{u}), \quad 
\text{and $\|u(T)-\uh\|$ small,}}
\end{subequations}
where $\hat{u}$ is a steady state of \reff{zivp}, and $E_s(\uh)$ 
is spanned by the eigenvectors of $-\pa_u G$ belonging 
to eigenvalues with negative real parts. As in 0D we then
need the SPP, i.e. 
\huga{\label{SPP} 
\dim E_s(\hat u)=\dim E_u(\hat u), 
} 
to chose an arbitrary initial point in the state space. For a  
discussion and possible extension of the SPP to PDEs see
App.~\ref{sec:spppdes}.

Arguably, the simplest discretization for \reff{mola}, at least in 1D, 
is a finite difference (FD) scheme, which has the advantage that we can directly
use \OCMAT{} for \reff{zivp},\reff{asmf-cond}, see \cite{grass2014}.  
However, here we opt
for a FEM ansatz for \reff{mola}, using the setting of \pdep\
\cite{p2pure, p2p2}
for two reasons: 
\bcen
\item We want to consider \reff{sldiff1} and hence \reff{slcan2} 
also on general 2D domains,
  and for more general models where the state variables may be
  vector valued functions already, see \S\ref{dsec}, again in 1D or 2D. In
  all these cases, FD and the coding of the respective spatially
  discrete systems may become rather inconvenient, while \pdep\
  provides convenient interfaces precisely for such systems.
  Moreover, for more complicated systems adaptive meshes may become
  important, which are more easily handled in a FEM discretization,
  and are already an integral part of \pdep.
\item As explained above, the canonical system may have many
  stationary states; it is thus desirable to use a continuation and
  bifurcation package to conveniently find CSS. The goal then is to
  ``seamlessly'' link the setting of \pdep{} for stationary problems
  with BVP solvers for \reff{zivp}.  
\ecen

  On the other hand, a drawback of spatial FEM discretizations is
  that the associated evolutionary problems have the natural form
  \huga{\label{fema} M\dot u(t)=-Ku(t)+MF(u(t))=:-G(u(t)), 
} where
  $u$ corresponds to the nodal values, $M, K\in
  \R^{2n\times 2n}$ are called the mass
  matrix and the stiffness matrix, respectively, and $F:\R^{2n}\ra\R^{2n}$ is the
  nonlinearity. Thus, $M$ and $K$ are large but sparse; the
``-'' signs in \reff{fema} comes from the convention that 
\pdep{} discretizes $-\Delta$ as $K$ (positive definite). 
  The occurrence of $M$ on the left hand side of \reff{fema} means
  that it is not of the form \reff{zivp}, and creates problems for the
  usage of standard BVP solvers.  Of course, $M$ is non--singular,
  and hence \reff{fema} can be rewritten as
  \begin{equation}
    \label{Minvode}
    \dot u=-M^{-1}Ku+F(u), 
  \end{equation}
  where for speed we can pre--calculate $M^{-1}K$. However, $M^{-1}$
  is no longer sparse, and already for intermediate $n$ ($n>1000$,
  say) this results in slow computations and in particular memory
  problems when using standard BVP solvers, which sort the Jacobian
  $-M^{-1}K+\pa_u F$ from each time--slice $t_0,\ldots,t_m$
  into a big Jacobian for the BVP problem. This can be alleviated by 
providing an approximate Jacobian
  $\tilde{J}=-A+\pa_u F$, where $A$ approximates $M^{-1}K$ via
  lumping, i.e.\ we drop entries  
from $M^{-1}K$ below a certain size $\del$. 
Of course, there is a trade off between accuracy of  $\tilde{J}$ and 
the number of its nonzeros. 

We mainly experimented with
  the \mlab{} solvers {\tt bvp4c} and similar, for instance the
  adaptions of {\tt bvp4c} already implemented in \OCMAT, and the
  solver TOM \cite{mazS2002, MT2009}\footnote{see also \httpTOM}. 
It turned out that TOM worked best in the ``lumped'' setting. Moreover, 
TOM was easy to modify in an ad hoc way to handle $M$ on 
the left hand side, and a new official release of TOM is scheduled 
that also uses $M$ \cite{mpriv}. 
 For speedup it is advisable to avoid numerical differentiation and
  hence to pass a Jacobian function {\tt J=fjac(t,u)} to TOM. 
This is generically very easy as \pdep\ in most cases provides 
a fast and easy way to assemble $J$. See \cite{U15p2} for 
implementation details. 

\subsection{The algorithm for the computation of a path 
to a CSS satisfying the SPP}
\label{sec:AlgorithmForTheComputationOfAPathSatisfyingSPP}
To calculate canonical paths 
from a given state $P_0$ that connect to some 
CSS $\hat{u}$ with the SPP we want to solve the two-point BVP 
\begin{subequations}
\label{eq:ocpathbvp}
\hual{
	M\dot u(t)&=-G(u(t)),\quad t\in(0,T)\label{dyn}\\
    P_i(0)&=P_{0,i},\ i=1,\ldots,n,\quad\text{($n$ left BC),}\label{lbc}\\
    \asymmat(u(T)-\hat u)&=0\in\R^n\quad\text{($n$ right BC), }
\label{rbc}
}
where $\Psi\in R^{n\times 2n}$ encodes the projection onto the unstable 
eigenspace, i.e.~$\asymmat (u-\hat u)=0$ for $u\in E_s(\hat{u})$, and where 
$T$ is the chosen truncation time. 
The calculation of $\Psi$ at startup, which for large $n$ turns out to be one 
of the bottlenecks of the algorithm,  also gives
  a lower bound for the time scale $T$ via $T\ge \frac 1{-\Re\mu_1}$,
  where $\mu_1$ is the eigenvalue with largest
  negative real part, i.e., gives the slowest direction of the stable
  eigenspace of $\hat{u}$. In our simulations we typically use
  $T$ between 50 and 100. 

In general, a BVP solver needs a good initial guess of $t\mapsto u(t)$ 
to solve problem \reff{eq:ocpathbvp}. Therefore we embed 
problem \reff{eq:ocpathbvp} into a family of problems replacing \reff{lbc} by
\begin{equation}
\label{eq:lbc_alt}
	P(0)=\al P_{0}+(1-\al)\hat P,\ \al\in[0,1],
\end{equation}
\end{subequations}
where we assume that for some $\al$ the solution is known: 
this holds for instance for $\al=0$ with the trivial 
solution $u\equiv \uh$. We may then gradually increase $\al$, 
using the last solution as the new initial guess. 
This is implemented in the algorithm summarized in Table \ref{tab2}. 

There are more sophisticated variants of the simple continuation in
Step 2 of Table \ref{tab2} (some of which are implemented in \OCMAT),
but the simple version in general works well for the problems we
considered. Nevertheless, it may be that no solution of
\reff{dyn}, \reff{rbc} and \reff{eq:lbc_alt} is found for $\al>\al_0$
for some $\al_0<1$, i.e., that the continuation to the intended
initial point fails. In that case usually the BVP problem undergoes a
fold bifurcation. We then use an adapted continuation step
2$^\prime$ that allows us to continue solutions around the fold. 

\begin{table}[ht]\begin{tabular}{c}
\hline\\[0mm]
\begin{minipage}{160mm}
      \bcen
    \item[Step 0 (selection of $\uhat$ and implementation of
      \reff{rbc}).]  Given $\uhat$  we solve the generalized adjoint 
EVP $\pa_u G(\uhat)^T\Phi=\Lambda M \Phi$ for the eigenvalues $\Lambda$ 
and (adjoint) eigenvectors $\Phi$, which also gives the defect $d(\uh)$. If 
$d(\uh)=0$, then from $\Phi\in \C^{2n\times 2n}$ 
we generate a real base 
of $E_u(\hat u)$ which we sort into the matrix  $\Psi\in \R^{n\times 2n}$.
    \item[Step 1 (selection of initial mesh and initial guess).]  
To start the BVP
      solver we choose the initial guess $u(t)\equiv \uh$ on a suitable initial
      grid $0=t_0<t_1<\ldots<t_m=T$. Typically, we choose $m=20$ at
      startup, and afterwards TOM uses its own
      mesh-adaption strategy. 
    \item[Step 2 (solution and continuation).]  Using \reff{eq:lbc_alt} 
we try to increase
      $\al$ in small steps $\delta$ to $\al=1$, in each step using the
      previous solution as the new initial guess, often trying 
$\delta=1/4$.  After thus having computed the first two solutions 
we may use a secant predictor for the subsequent steps. 
    \item[Step 2$^\prime$ (arclength continuation).] If the 
continuation fails for $\alpha>\alpha_0$ with
      $\alpha_0<1$, then we use a pseudo--arclength continuation for 
a modified BVP, letting $\alpha$ be a free
      parameter. \ecen 
\text{}\\[-4mm]
    \end{minipage}\\\hline
\end{tabular}
  \caption{The continuation--algorithm {\tt iscont} 
(Initial State Continuation);  Steps 0,1 are preparatory, Step 2 
    or 2$^\prime$ is repeated. See \cite{U15p2} for implementation details, 
and for remarks on performance. 
\label{tab2}}
\end{table}

\brem\label{tnormrem}{\rm (a) For some applications it is useful to 
rescale the time $t=T\tau$ and hence consider $M\dot u(\tau)=TG(u(\tau)))$ 
on the normalized time interval $\tau\in[0,1]$, which turns the truncation 
time $T$ into a free parameter. This is for instance 
implemented in {\tt OCMat}, but for simplicity not used here.

(b) Similarly to the normalized normalized objective value $J_{ca}$ in 
(\ref{sldiff1}b), in the bifurcation diagrams we use the normalized 
$L^2$ norm for comparison between different domains and space dimensions, 
i.e., henceforth, $\|P\|_2:=\|P\|_{L^2}/\sqrt{|\Om|}$, 
and in the table in Fig.~\ref{f1} we present averaged values, i.e., 
\huga{\label{avdef} 
\spr{P}:=\frac 1 {|\Om|} \int_\Om P(x)\dd x, \quad 
\spr{k}:=\frac 1 {|\Om|} \int_\Om k(x)\dd x. 
}
To take the finite truncation time $T$ into account we let 
\begin{equation}
\label{eq:numobjval}
	\tilde J(k(\cdot),T):=J(P_0(\cdot),k(\cdot,\cdot))
+\frac{\E^{-r T}}{r}J_{ca}(P(T),k(T)). 
\end{equation}
Obviously, for $T\gg\dfrac{1}{r}$ the last term can be made arbitrarily 
small, while for CSS it yields the exact (discounted) objective value. 
In the following we drop the tilde in \reff{eq:numobjval}, 
and write, e.g., $J_{\hat u}$ for the objective value of 
a \css\ $\hat u$, and, e.g., 
$J_{P_1\to \hat u}$ for the CS--path which goes from 
$P_1$ to $\hat u$. 
\eex}\erem 

\subsection{1D canonical steady states}\label{1dsec}
Recall that first we use \PDETOPATH{} to study the steady state problem 
for \reff{zivp}, i.e.,
\begin{subequations}
  \label{slepde}
  \begin{align}
    &0=-\frac{1}{q(x)}-bP(x)+\frac{P(x)^2}{1+P(x)^2}+
    D\Delta P(x),\label{sldiffusion_model_cansys_epde1}\\
    &0=2cP(x)+q(x)\left(r+b-\frac{2P(x)}{\left(1+P(x)^2\right)^2}\right)
    -D\Delta q(x),\label{sldiffusion_model_cansys_epde2}\\
    &\evalat{\pa_\nu
      P(x)}_{\pa\Omega}=0,\quad\evalat{\pa_\nu q(x)}_{\pa\Omega}=0.
    \label{sldiffusion_model_cansys_epde3}
  \end{align}
\end{subequations}
In 1D we choose $\Om=(-L,L)\times(-\del_y,\del_y)$ with 
Neumann BC on all boundaries and small $\del_y$ such that we can  
use just two grid--points in $y=x_2$--direction and the solutions will
be constant in $y$; this we call a quasi 1D setup.  
The steady states of the canonical system for the 0D
model \labelcref{sl_model} are \fcss{} of \crefsys{slepde}, and for easy 
reference we introduce the acronyms in Table \ref{tab1b}. 
\begin{table*}[h]
\bce 
\begin{tabular}{p{10mm}p{140mm}}
name&description\\
\hline
FSM&Flat State Muddy, the upper FCSS branch with a high phosphor load $P$\\
FSI&Flat State Intermediate, the upper half of the second FCSS branch, intermediate $P$\\
FSC&Flat State Clean, the lower half of the second FCSS branch, low $P$\\
\hline 
\end{tabular}
\caption{Classification of the FCSS branches, see also 
Fig.\ref{f1}. The high $P$ state is also called eutrophic, and our ``muddy'' 
refers to the fact that under eutrophic conditions there are a lot of algae 
and other organic matter in the lake, while under oligotrophic 
conditions (low $P$) 
the water is much ``cleaner''. \label{tab1b}}
\ece 
\end{table*}
The FCSS are of course independent of the domain, but to search 
for \hcss{} bifurcating from \fcss{}, 
the domain size $2L$ should be close to a multiple of $2\pi/k_c$,
where $k_c$ is the wave number of a Turing bifurcation. The parameters
in \reff{parsel}, with $b$ near 0.7, yield $k_c\approx 0.44$
\cite{BX08}, and for simplicity we then choose $L=2\pi/0.44\approx14.28$. 

\begin{figure}[ht]
\begin{tabular}{p{0.31\textwidth}p{0.3\textwidth}l}
(a) BD of CSS, (normalized) $L^2$ norm over $b$&(b) BD, current values $J_{c,a}$& (c) example CSS\\
\ig[width=0.3\textwidth]{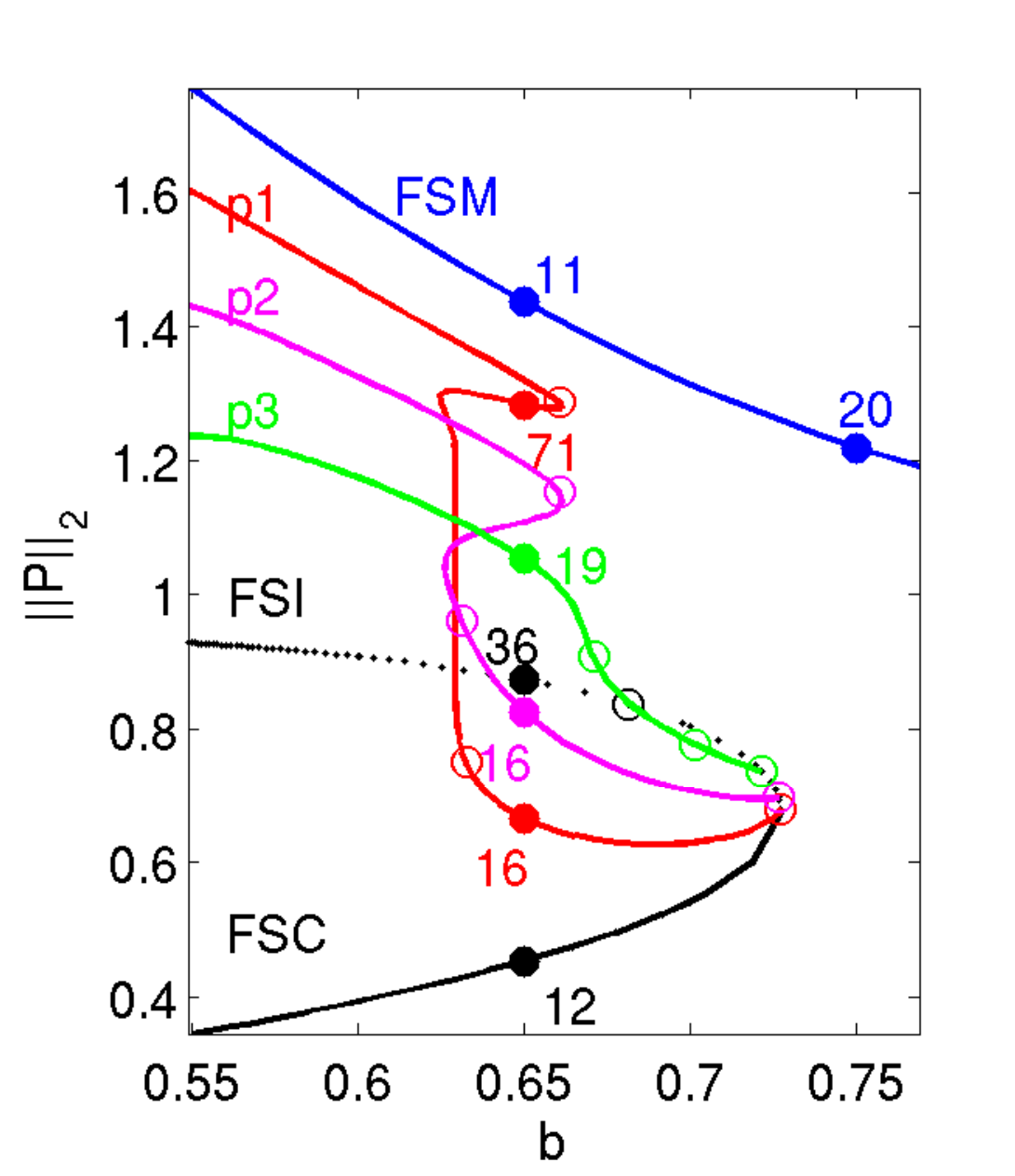}&
\ig[width=0.3\textwidth,height=60mm]{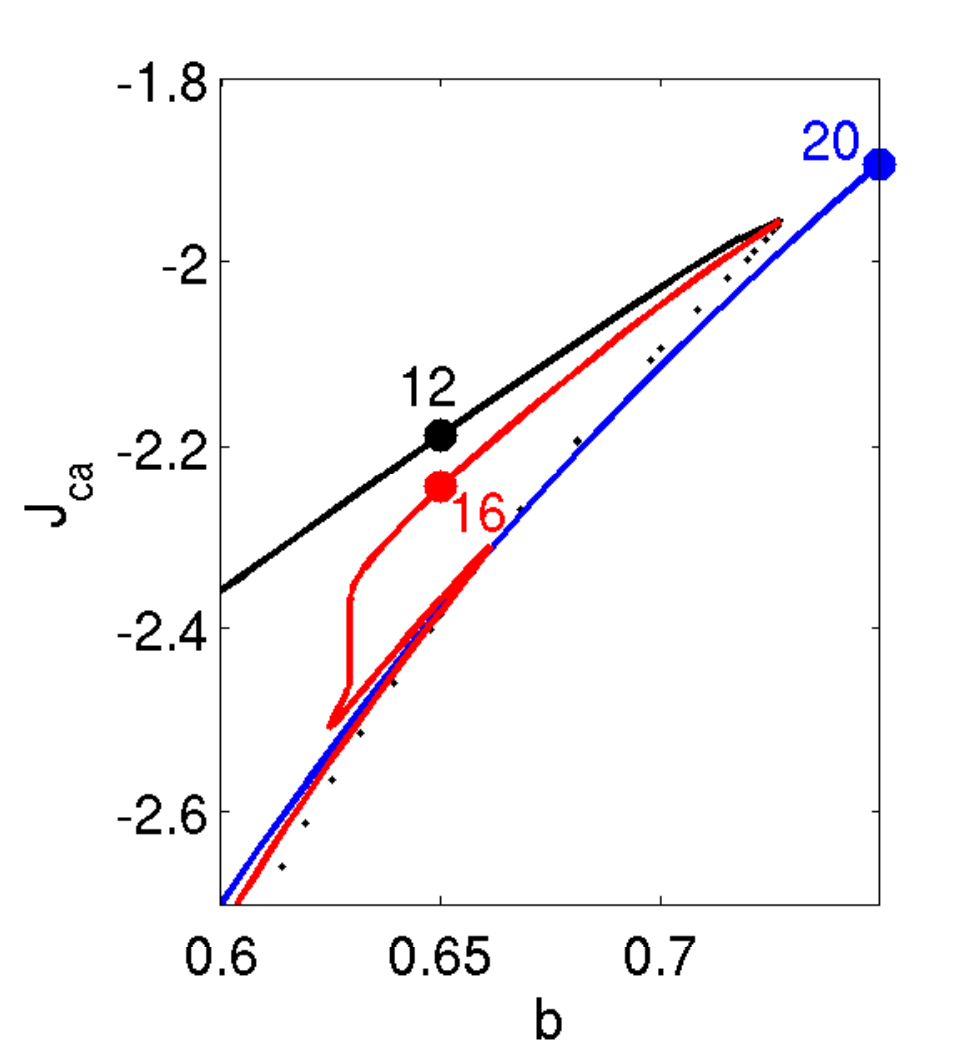}&
\raisebox{30mm}{\begin{tabular}{l}
\ig[width=0.15\textwidth]{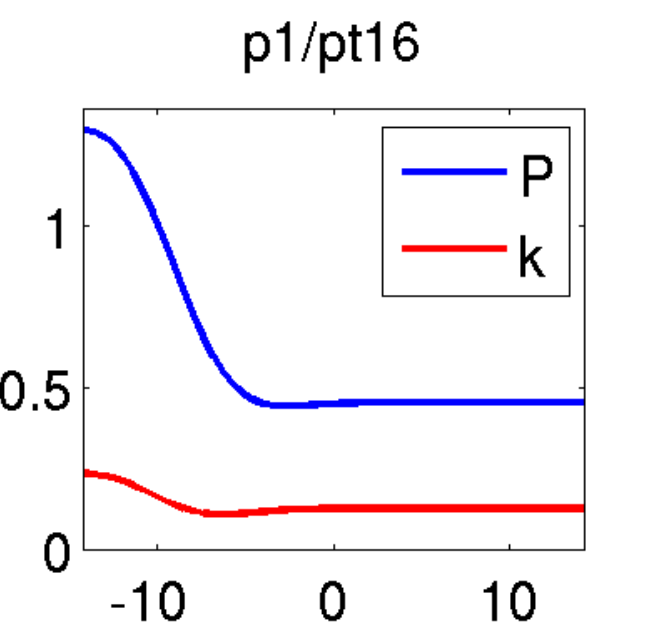}\ig[width=0.15\textwidth]{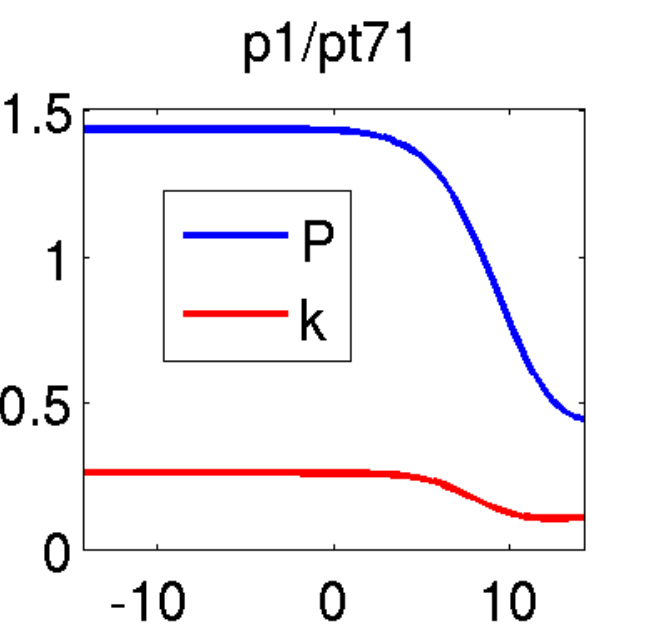}\\
\ig[width=0.15\textwidth]{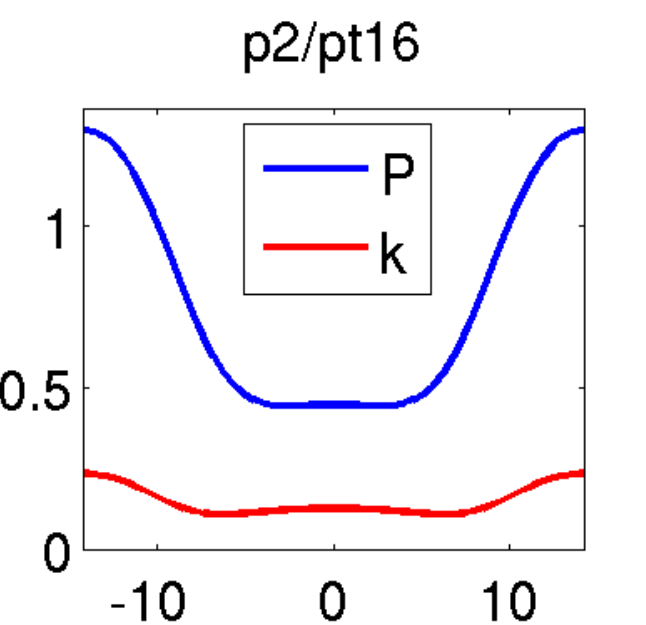}\ig[width=0.15\textwidth]{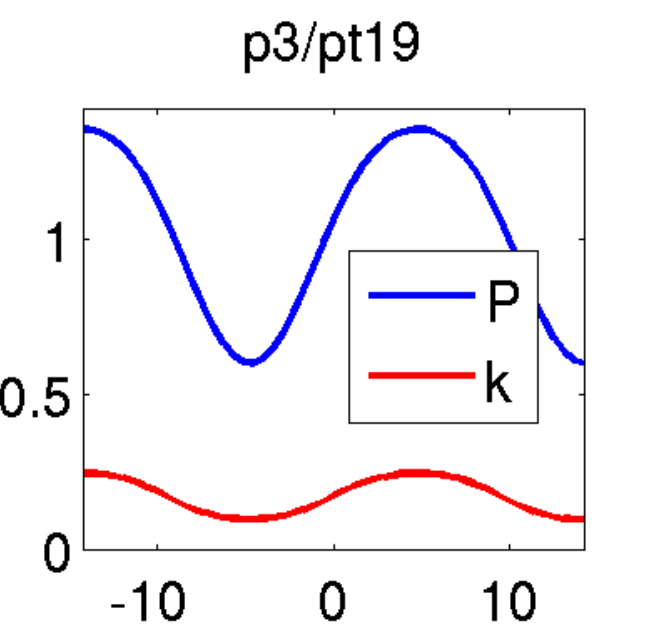}
\end{tabular}}
\end{tabular}

\bce 
(d) Characteristics of points in (a)-(c).\\
\begin{tabular}{p{20mm}cccc|p{20mm}cccc}
name&$\spr{P}$&$\spr{k}$&$J$&$d$& name&$\spr{P}$&$\spr{k}$&$J$&$d$ \\
\hline
FSM/pt20& 1.22 & 0.32 & -63.11 & 0&p1/pt16& 0.61 & 0.14 & -74.83 &-1\\
FSM/pt11& 1.44 & 0.26 & -79.28 & 0&p1/pt71& 1.24 & 0.22 & -78.93& 0\\
FSI/pt36&0.87 & 0.13 & -79.47&-5 & p2/pt16&0.76 & 0.15 & -76.70 & -2\\
FSC/pt12&0.45 & 0.12 & -72.95 &0&p3/pt19&1.02 & 0.17 & -79.48 &-3\\ 
\hline 
\end{tabular}
\ece 
   \caption{Basic bifurcation diagrams (a) and (b), example plots (c), and 
characteristic values of selected CSS (d). In (a), the blue and black lines 
represent the FCSS, and for instance the red line {\tt p1} contains patterned 
CSS with one ``interface'' between high and low $P$. The numbered points 
on all these lines correspond to selected solutions plotted in (c), and 
characterized in (d). The small circles in (a) denote bifurcation points. 
The values $J_{ca}$ and $J=J_{ca}/r$ of the CSS 
are all negative, but this is merely 
a question of offset. } 
  \label{f1}
\end{figure}

Continuing the FSI branch in $b$ we find a number of Turing like 
bifurcations to PCSS, and follow four of these; see Fig.~\ref{f1}a,b for the
bifurcation diagram(s), and (c) for example
solutions\footnote{The notation, e.g., {\tt p1/pt16} follows the 
\PDETOPATH{} scheme, e.g.: continuation step 16 
on the branch p1 is stored in folder {\tt p1}
  and file {\tt pt16.mat}.}. On the branches p1,p2, and p3 there are secondary 
bifurcations, not further considered here. 
 For the subsequent examples we focus on $b=0.75$ and $b=0.65$, 
check the SPP \reff{SPP} for all CSS 
and find it only to be fulfilled for the FSM, for the FSC, and
for some points on the p1 branch, e.g.~at point 71 
between the folds (see Fig.~\ref{f1}).

\subsection{1D canonical paths} 
For $b=0.75$, the only CSS is the FSM ({\tt FSM/pt20}). It has the 
SPP, we can reach it from an arbitrary initial state $P_0$, 
and thus it is a globally stable \foss. Therefore, this regime is not very 
interesting, and we immediately turn to the case $b=0.65$ with 
multiple CSS.

For $b=0.65$ seven CSS are marked in Fig.~\ref{f1}a, and 
characterized in the table, where only three satisfy
the SPP. These are the FSC (which has the maximal value among these CSS), 
the FSM, and the \hcss{} {\tt p1/pt71}, subsequently denoted as 
$\uh_{{\rm PS}}$. 
Next we numerically analyze which of these \css\ belong to optimal paths. 

\subsubsection{\hcss{} not satisfying the SPP}
\label{nonSPPsec}
From the analysis of non-distributed optimal control problems
we know that steady states that do not satisfy the SPP can
nevertheless be optimal. To illustrate that this is at least not typical 
in the SLOC model, we  
first compare the objective value of some (arbitrary chosen) 
CSS, e.g., the \hcss{}
$\uh_{{\rm PS}^-}(\cdot)$ {\tt p3/pt19}, not satisfying
the SPP, with that of the $t$--dependent canonical paths 
$ u_i(\cdot,\cdot)$  which 
connect to $\hat u_i(\cdot),\ i\in\{{\rm
    FSC}, {\rm FSM}, {\rm PS}\}$. 
For the objective values we write 
$J_{{\rm PS}^-}$ for the \css, and, e.g., 
$J_{{\rm PS}^-\to{\rm FSC}}$ for the canonical path which goes from 
$P_{{\rm PS}^-}$ to $\uh_{{\rm FSC}}$. 
The optimal solution for $P_0=\phossnspp$ has to
satisfy
\begin{equation}\label{Vord}
  V(\phossnspp)=\max\{J_{{\rm PS}^-},\ J_{\phossnspp\to{\rm FSC}},\ 
J_{\phossnspp\to{\rm FSM}},\ J_{\phossnspp\to{\rm PS}}\}. 
\end{equation}
The canonical paths are given in Fig.~\ref{ocf3}a--c, while (d) presents 
some norms along the path in (a), 
which show that and how fast $u(t)$ (including the co--states $q$) converges 
to $\uh$. In all cases we find without problems canonical paths to 
both \fcss{} and the \hcss; 
in particular the path to FSM is rather quick. 

For the objective values we find, up to 2 significant digits, 
\hual{\label{jord}
  J_{{\rm PS}^-}{=}-79.48<J_{{\rm PS}^-\to{\rm FSC}}{=}-78.24
<J_{{\rm PS}^-\to{\rm PS}}{=}-78.19<J_{{\rm PS}^-\to{\rm FSM}}{=}-77.5 . 
}
Thus, the optimal path is the one converging to $\uh_{{\rm FSM}}$. 
Repeating these steps for every \hcss{} not satisfying the SPP we find that
these are always dominated by paths converging to one
of the \fcss. Therefore, only $\uh_{{\rm FSC}}$, $\uh_{{\rm FSM}}$ and 
$\uh_{{\rm PS}}$ remain as candidates for \oss. 

{\red Before we turn to determining OSS, we briefly discuss 
the canonical paths in Fig.~\ref{ocf3}. We focus on (a), but similar 
remarks apply to (b,c) as well. In (a), the initial $P(\cdot,0)$ is above the 
target $\hat P_{\text{FSC}}$, so naively we may expect that the control 
$k$ should start below the target $\hat k_{\text{FSC}}$ and slowly increase to 
$\hat k_{\text{FSC}}$. However, such a control would not be optimal. 
Instead, $k$ is initially similar to $\hat k_{\text{PS}^-}$, and 
in particular $k(\cdot,0)$ is large where $P_0(\cdot)$ is already large. 
Only after a short transient $k$ drops below $\hat k_{\text{FSC}}$ and 
then behaves as expected. 

\begin{figure}[ht]
\begin{tabular}{ll}
(a) canonical path from {\tt p3/pt19} to FSC
&(b) canonical path from {\tt p3/pt19} to FSM\\
\ig[width=0.245\textwidth]{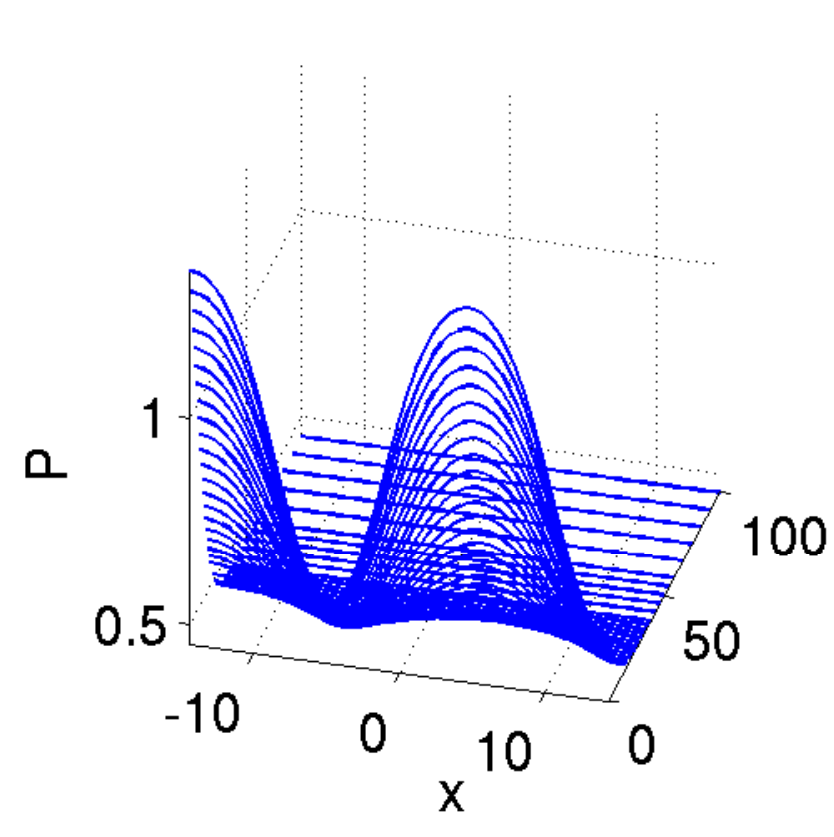}
\ig[width=0.245\textwidth]{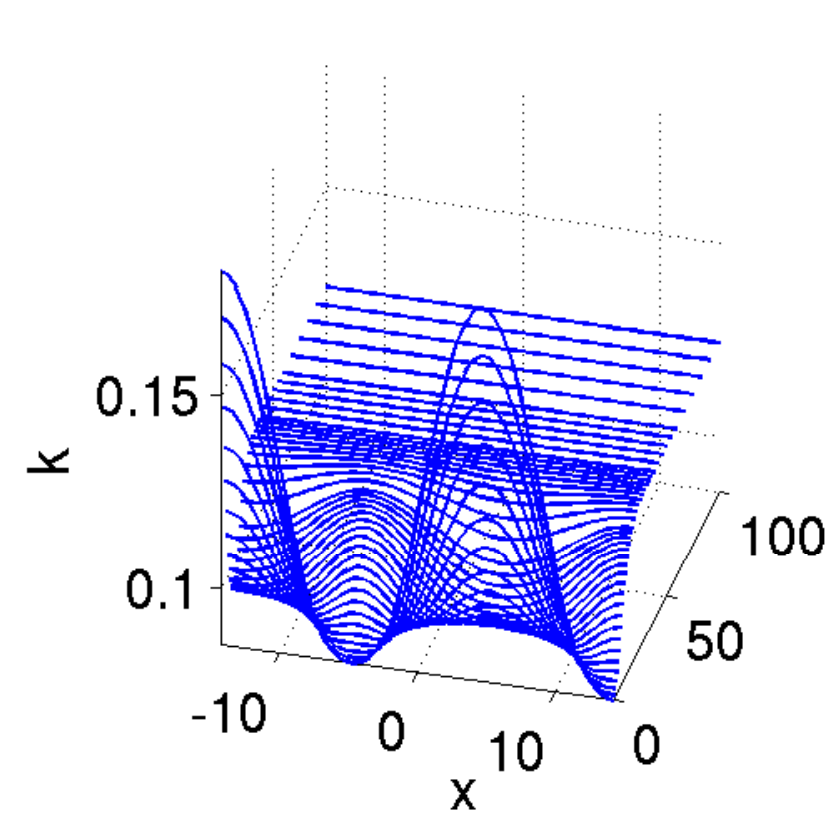}&
\ig[width=0.245\textwidth]{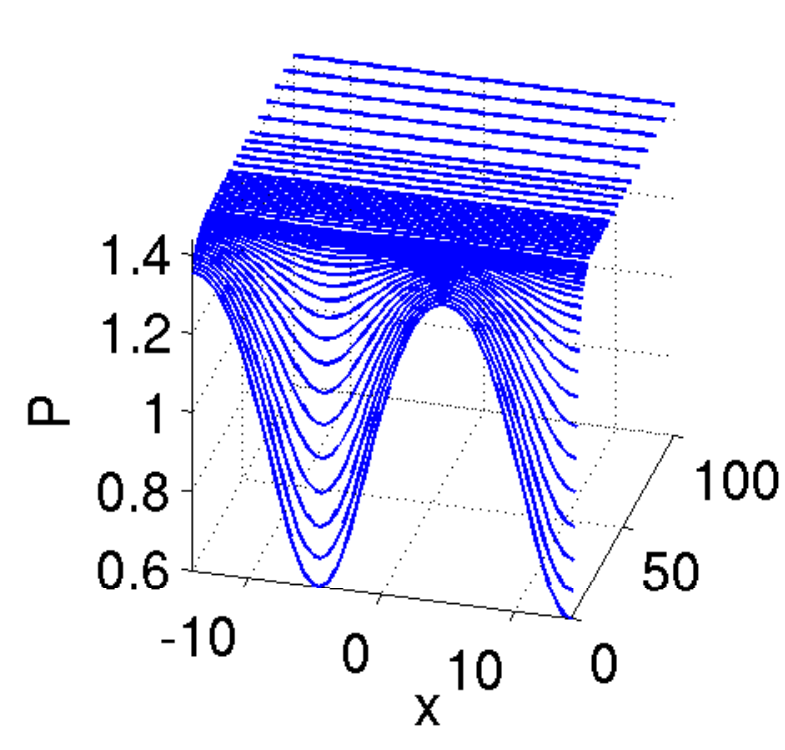}
\ig[width=0.245\textwidth]{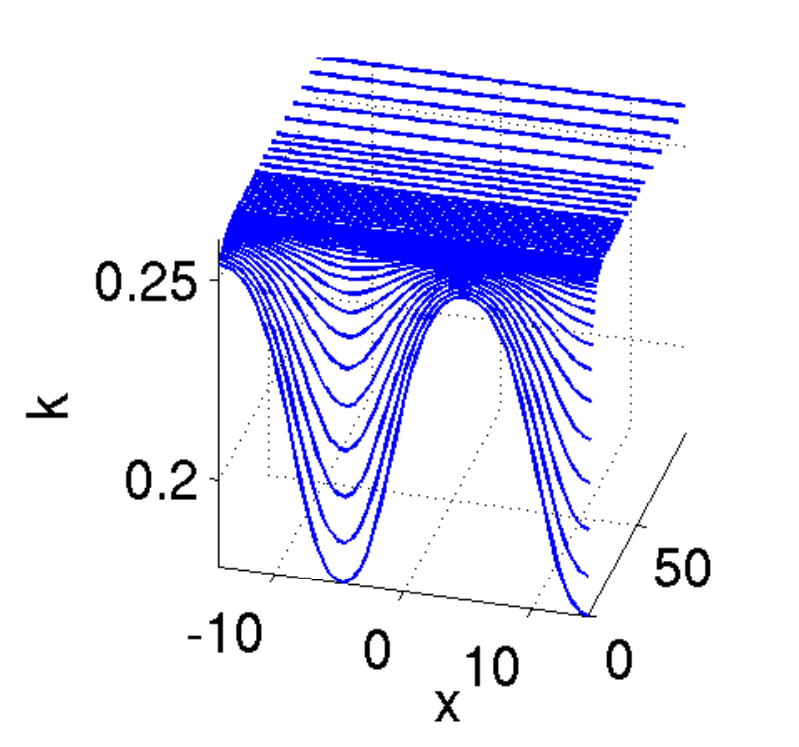}\\
(c) canonical path from {\tt p3/pt19} to PS&
(d) diagnostics for a)\quad (e) $P$ for a naive control \\
\ig[width=0.245\textwidth]{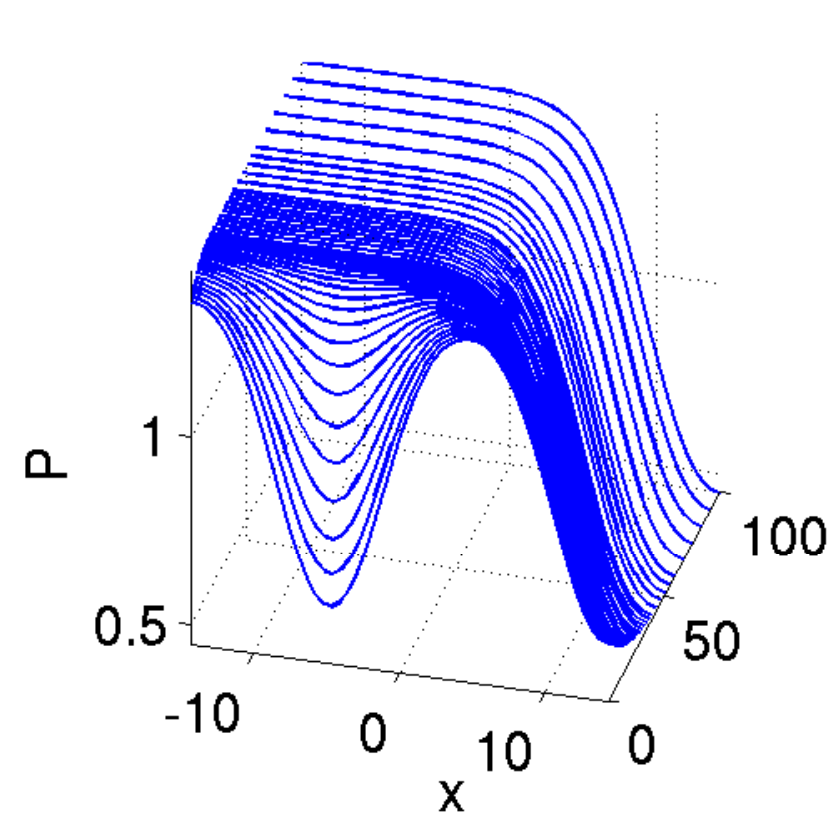}
\ig[width=0.245\textwidth]{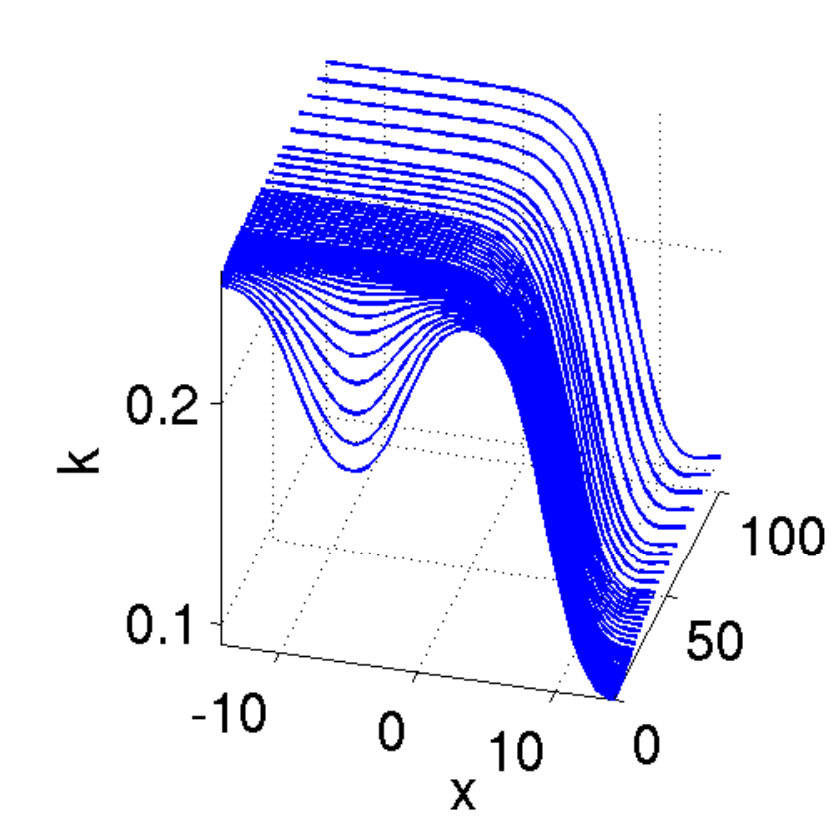}&
\ig[width=0.23\textwidth]{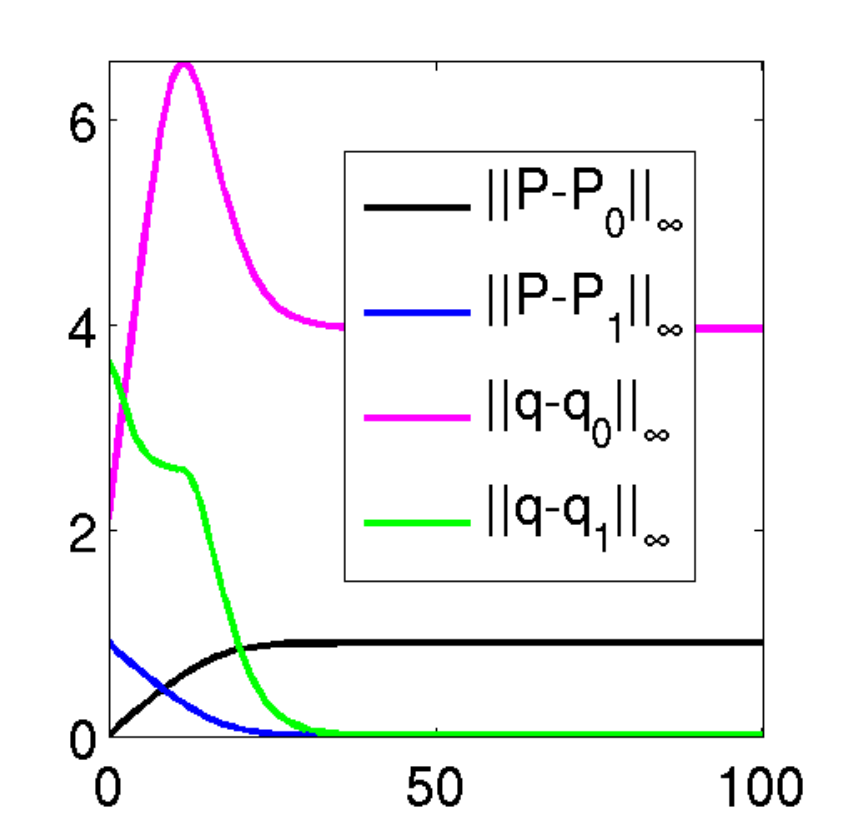}
\ig[width=0.245\textwidth]{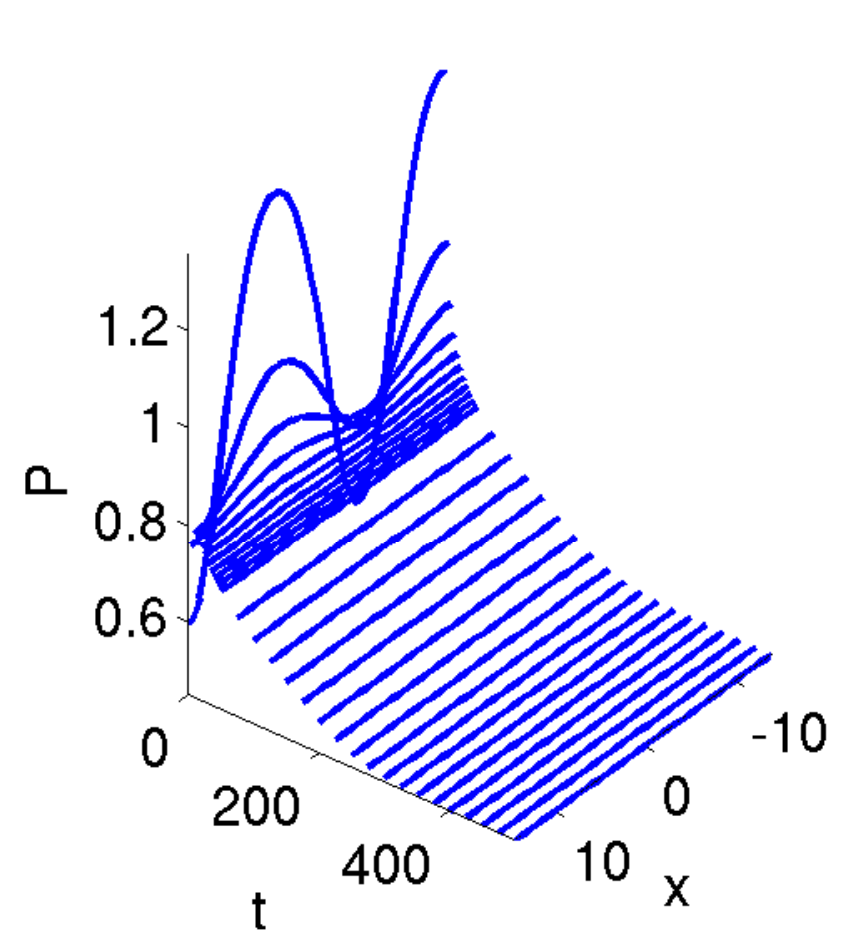}
\end{tabular}
   \caption[]{Canonical paths from the (state values of) 
\hcssnspp{} {\tt p3/pt19} to 
FSC (a), FSM (b), and PS (c) at $b=0.65$, and typical path diagnostics (d). 
For comparison, (e) shows the solution $P$ of the initial value problem 
(\ref{sldiffusion_model_dyn}) with $P(0)$ as in (a) and the externally chosen control 
$k(t)\equiv k_{\text{FSC}}$ for all $t$. }
  \label{ocf3}
\end{figure}

At first sight, this startup behavior of $k$ may appear rather 
counter--intuitive. However, the reason is that we do {\em not} want to 
drive the system to $\hat u_{\text{FSC}}$ ``as quickly as possible'', 
which 
essentially would amount to choosing $k$ ``as small as possible'' at startup. 
Instead, we want to maximize $J$, and for this it pays off to have, 
for a short transient, $k$ large near the maxima of $P(\cdot,0)$. 

To illustrate this point, in Fig.~\ref{ocf3}(e) we {\em choose} the naive control 
$k(t)\equiv k_{\text{FSC}}$ for all $t$ and numerically integrate 
the {\em initial value problem} (\ref{sldiffusion_model_dyn}). 
While this does take us to $\uh_{\text{FSC}}$, the first observation is that 
this needs a rather long time. Secondly, for the value of this 
solution we obtain $J=-80.1$, which is even worse than the starting CSS 
with $J_{{\rm PS}^-}{=}-79.48$. }

\subsubsection{Determining optimal steady states} 
\label{OCPsec} 
We return to the question whether one or more of the CSS at $b=0.65$ 
 with the SPP from \S\ref{nonSPPsec} are optimal, and 
proceed in three steps. First we search for a canonical path
starting at $\phoss$ and connecting to
$\uh_{{\rm FSC}}$. In the second step we repeat
that for $\uh_{{\rm FSM}}$, and in the last step we
check if both or only one of the \fcss{} are 
optimal.\footnote{The second step reveals that the first
  step is superfluous, but this we do not know a priori.}

\begin{figure}[ht]
  \centering \subfloat[$J_{i,a}$ for FSC target]
{\ig[width=0.23\textwidth]{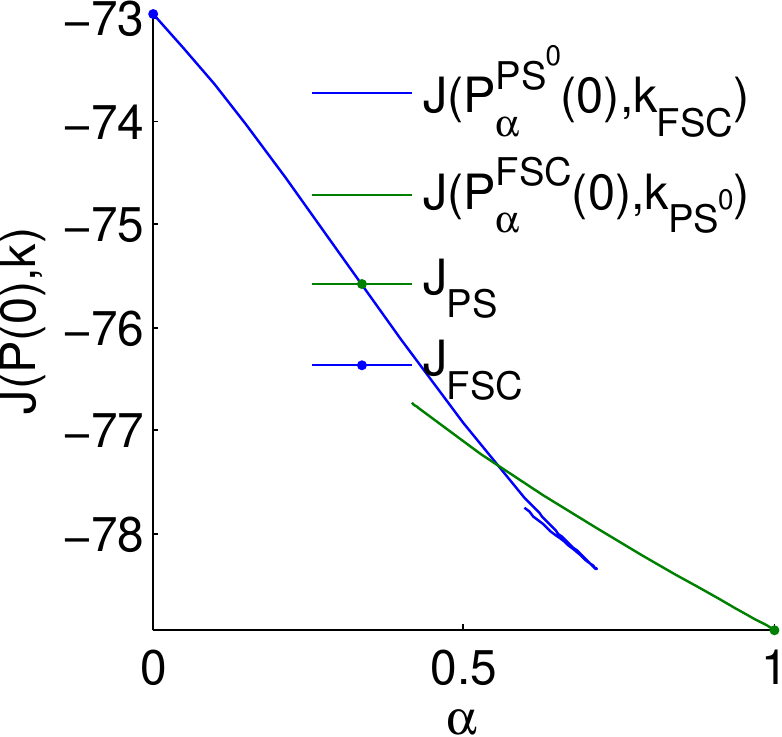}
\label{fig:skibacases1}}\ 
\subfloat[Two paths to FSC for $\al=0.6$]
{\ig[width=0.18\textwidth,height=30mm]{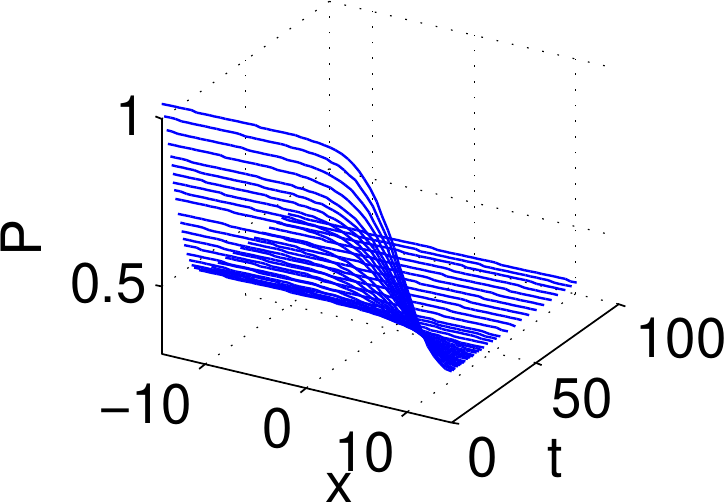}
\ig[width=0.18\textwidth,height=30mm]{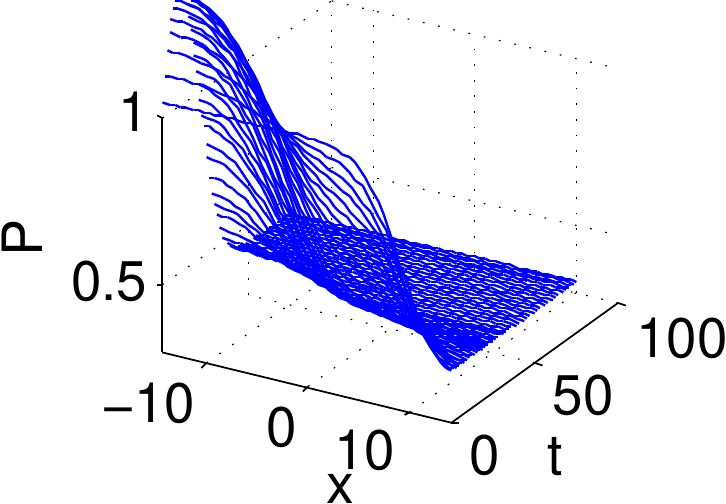}
\label{fig:skibacases3}}\ 
\subfloat[Diagn.~for left and right paths in (b)]
{\ig[width=0.17\textwidth,height=35mm]{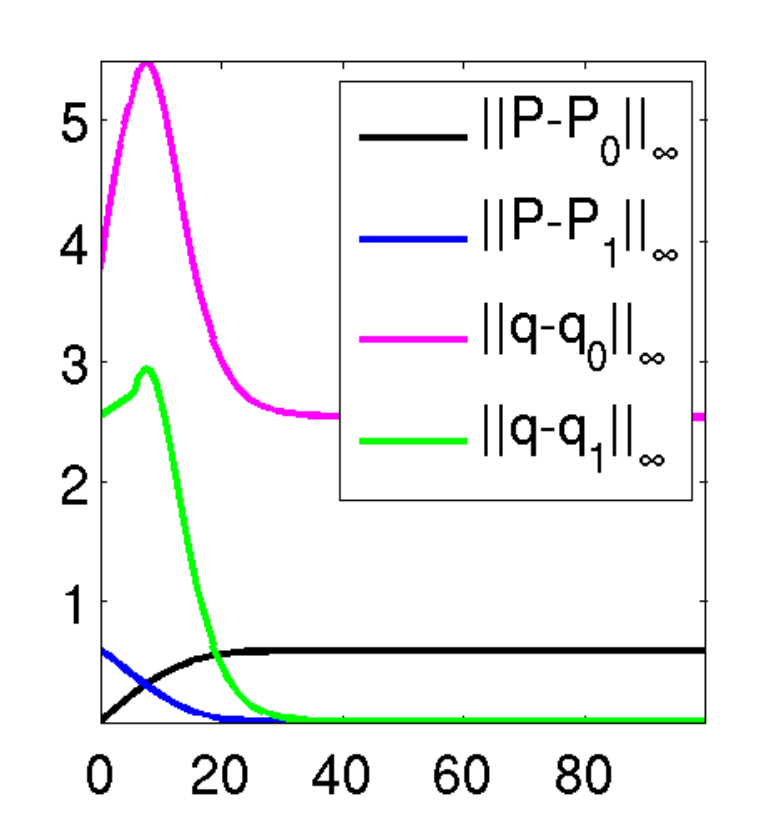}
\ig[width=0.17\textwidth,height=35mm]{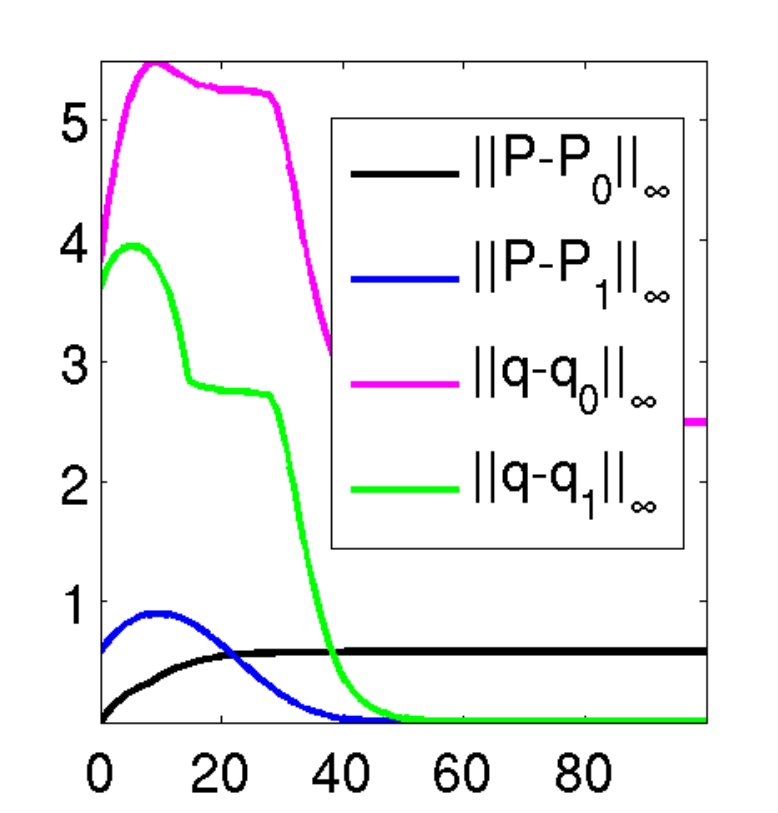}
\label{fig:skibacases3b}}\\
\subfloat[A Skiba candidate between $P_{{\rm PS}}$ and FSC]
{\ig[width=0.18\textwidth,height=30mm]{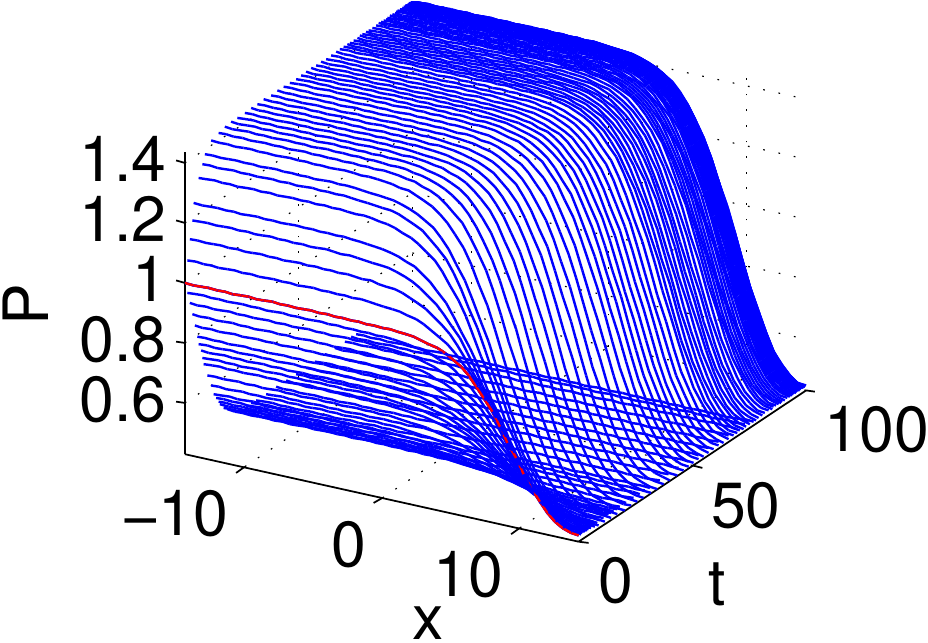}
\ig[width=0.18\textwidth,height=30mm]{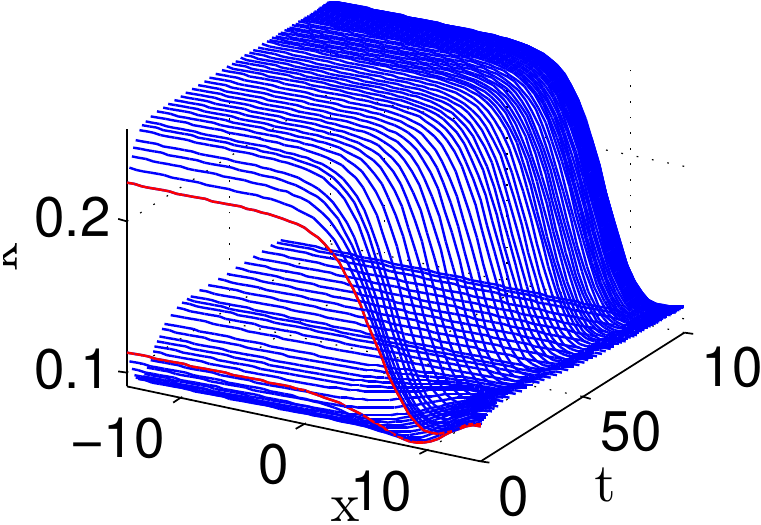}
\label{fig:skibacases4}}
  \subfloat[$J_{i,a}$ for FSM target]
{\ig[width=0.23\textwidth]{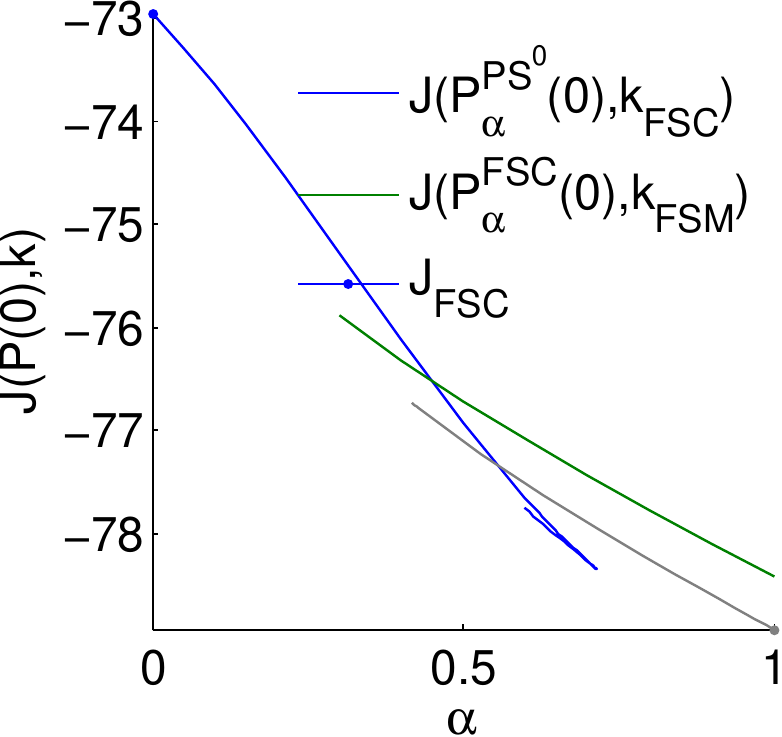}\label{fig:skibacases21}}\quad 
\subfloat[Dominating path to FSM]
{\ig[width=0.24\textwidth,height=30mm]{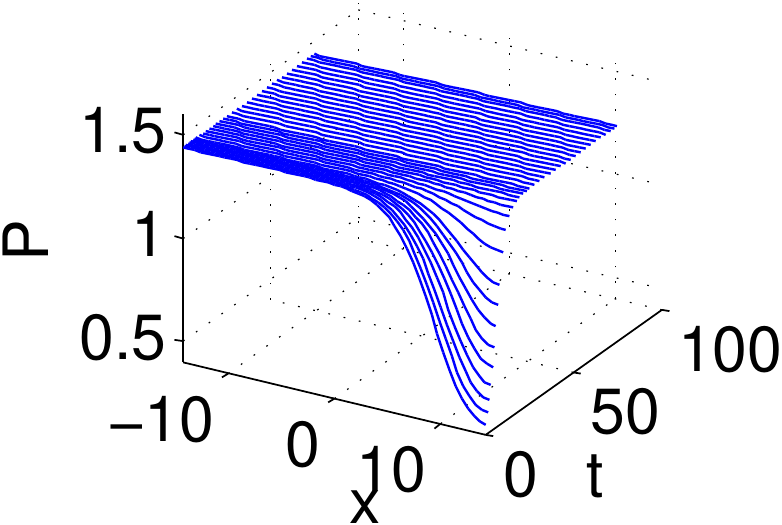}
\label{fig:skibacases2}}\\
 \subfloat[0D Shallow Lake]
{\ig[width=0.23\textwidth]{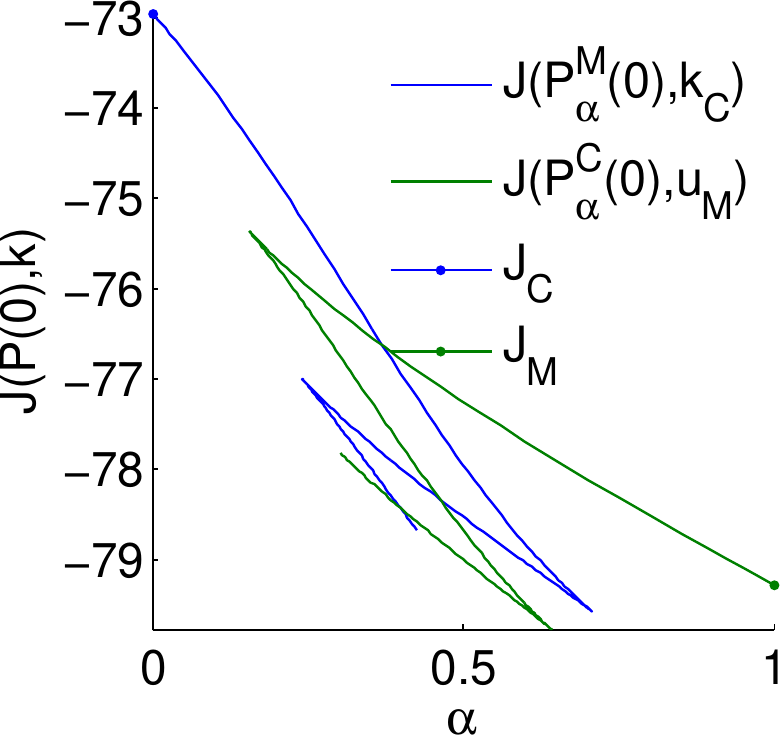}
\label{fig:skibacases30}}
  \subfloat[Homog.~Skiba point and paths to FSC and FSM]
{\ig[width=0.225\textwidth]{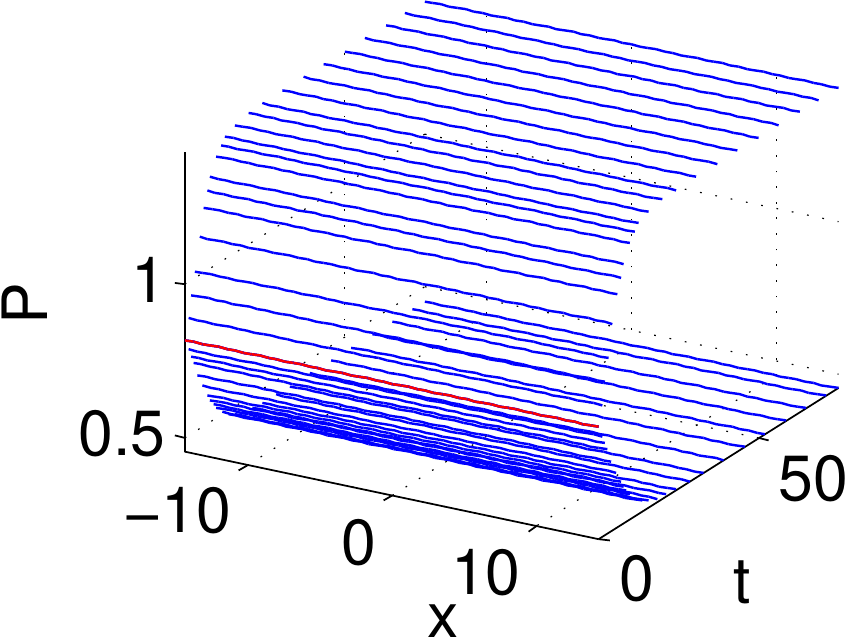}
\label{fig:skibacases6}}\quad 
  \subfloat[A patterned Skiba point between FSC and FSM]
{\ig[width=0.225\textwidth]{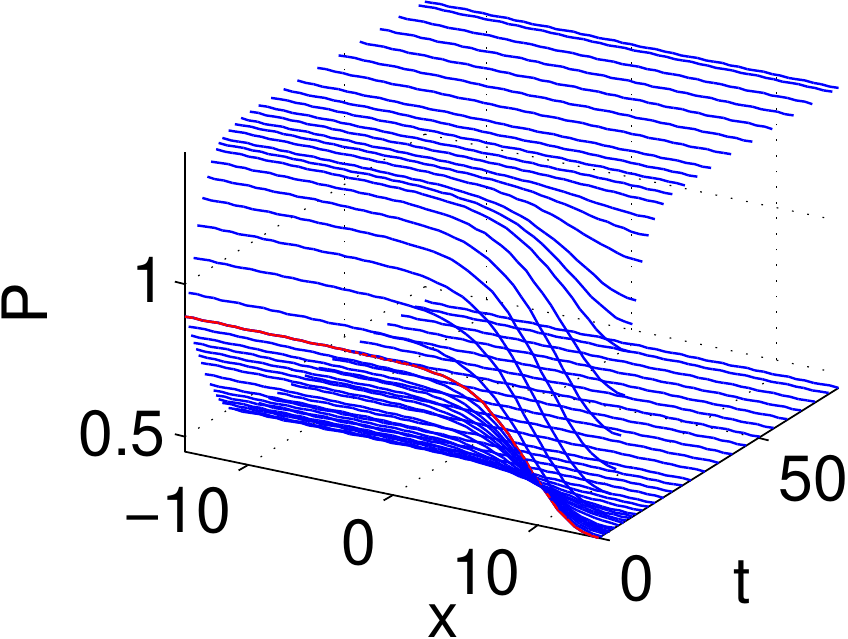}
\ig[width=0.225\textwidth]{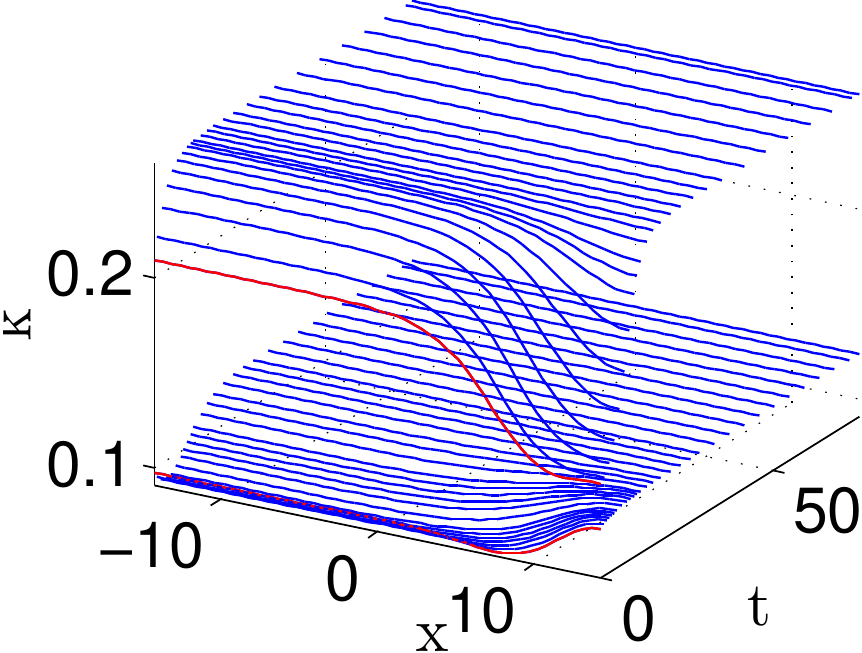}
\label{fig:skibacases5}}
  \caption{Canonical paths to various CSS for $b=0.65$, and illustration of 
some Skiba points; see text for details. 
  \label{fig:skibacases}}
\end{figure}

\paragraph{Paths between $\phoss$ and $\uh_{\text{FSC}}$ -- 
a Skiba candidate.}
\label{sec:FOSSPfossoneAndHOSSPhoss}
Using {\tt iscont} to get a canonical path starting at $\phoss$ and
converging to the FSC it turns
out that the continuation \reff{eq:lbc_alt},  
i.e., 
\huga{\label{al17}
P_\al (0):=\al P_{{\rm PS}}+(1-\al)P_{{\rm FSC}}, 
}
yields a fold around $\al\approx 0.6$, see the blue
curve in Fig.~\ref{fig:skibacases1}), and that no canonical path 
connecting $\phoss$ to $\uh_{{\rm FSC}}$ exists. Instead, multiple
solutions that converge to $\uh_{{\rm FSC}}$ 
exist for initial distributions of the form \reff{eq:lbc_alt} 
with $\alpha\in[0.6,0.71]$; two examples 
are shown in Fig.~\ref{fig:skibacases3}, and their diagnostics in (c). 

Similarly, trying to continue to a path that 
connects $\pfossone$ and $\uh_{{\rm PS}}$ yields 
a fold (green curve in (a)), and no such path exists. 
However, the solutions returned
during the continuation process allow us to determine and compare the
respective objective values. This yields that there exists a specific
initial distribution where the objective values are equal, 
given by the intersection of the green and blue curves in 
Fig.~\ref{fig:skibacases1}. Thus, from an economic point of view both 
solutions are equal. This suggests that 
$P_S$ is a Skiba or indifference threshold point 
(distribution)\footnote{however, taking into account also the FSM solution, 
below we identify $P_S$ as only a Skiba candidate}, 
well known from non-distributed optimal control problems, see, e.g., 
\cite{wagener2003,kiselevawagener2008a}. The
paths $u$ starting at $P_S$ (red curve)
are depicted in Fig.~\ref{fig:skibacases4}: $P(\cdot,0)$ is the same 
 for both solutions, but the controls $k(\cdot,0)$ 
are different.   In any case, to assure that these solutions are optimal 
we have to prove that no other dominating solution
exists. Thus in a last step we calculate the objective values of 
the paths converging to $\uh_{{\rm FSM}}$.  

\paragraph{From $\phoss$ to $\uh_{{\rm FSM}}$.} 
\label{sec:FOSSPfosstwoAndHOSSPhoss}
Here the continuation is successful and we find a
path connecting  $\phoss$ to $\uh_{\text{FSM}}$. Comparing the objective values
 reveals that the 
\hcss{} is dominated by the solution converging to the FSM, 
see Fig.~\ref{fig:skibacases21} and \ref{fig:skibacases2}. Thus,
the \hcss{} is ruled out as an optimal steady state, and therefore 
we a posteriori identify $P_S$ as only 
a Skiba \emph{candidate} as it does \emph{not} 
separate two optimal steady states. 

\paragraph{A Skiba manifold between 
$\pfosstwo$ and $\pfossone$.}
It is well known that in 0D the FSC and FSM are only 
locally  stable with regions of attractions separated by a Skiba
manifold (parametrized by, e.g., $b$) of homogeneous 
solutions, \cite{kiselevawagener2008a}, see Fig.~\ref{fig:skibacases30} 
for our case $b=0.65$. Of course, this also yields a homogeneous 
Skiba distribution in 1D, see Fig.~\ref{fig:skibacases6}. More generally, 
we may expect the domains of attraction of the FSC and the FSM  
to be separated by an (in the continuum limit infinite dimensional) Skiba 
manifold $M_S$, for fixed $b$. 

A continuation process, analogous to the
  non-distributed case, see \cite{grass2012}, could be used to
  approximate this manifold $M_S$. However, 
to find a non homogeneous example point on that manifold, here we can 
readily combine 
Fig.~\ref{fig:skibacases1} and \ref{fig:skibacases21}, to find 
a Skiba distribution of the form 
\huga{
P_{{\rm Skiba}}=\al\hat P_{{\rm FSC}}+(1-\al)\hat P_{{\rm PS}}; 
}
see Fig.~\ref{fig:skibacases5} for paths to the FSC and the FSM yielding 
the same $J=-76.3$.

\paragraph{Summary for 1D.} The picture that emerges 
is as follows: for $b>b_{{\rm fold}}\approx 0.727$ 
the FSM as the only CSS is the globally stable FOSS, while 
for $b<b_{{\rm fold}}$ there exist multiple CSS. Specifically for 
$b=0.65$, $\uh_{{\rm FSC}}, \uh_{{\rm FSM}}$ 
and one of the PCSS have the 
saddle point property, but only $\uh_{{\rm FSC}}$ and 
$\uh_{{\rm FSM}}$  
are optimal, and in particular no POSS exists. The FOSS FSC and FSM are 
separated by a (presumably rather complicated) Skiba manifold $M_S$, 
and Fig.~\ref{fig:skibacases6} and  \ref{fig:skibacases5} show just two 
examples of points on $M_S$.

\subsection{Outlook: 2D results}\label{2dsec}
\begin{figure}[!ht]
{\small 
\begin{tabular}{p{45mm}p{90mm}}
(a) Two 2D PCSS branches&(b) some example plots\\
\ig[width=45mm,height=44mm]{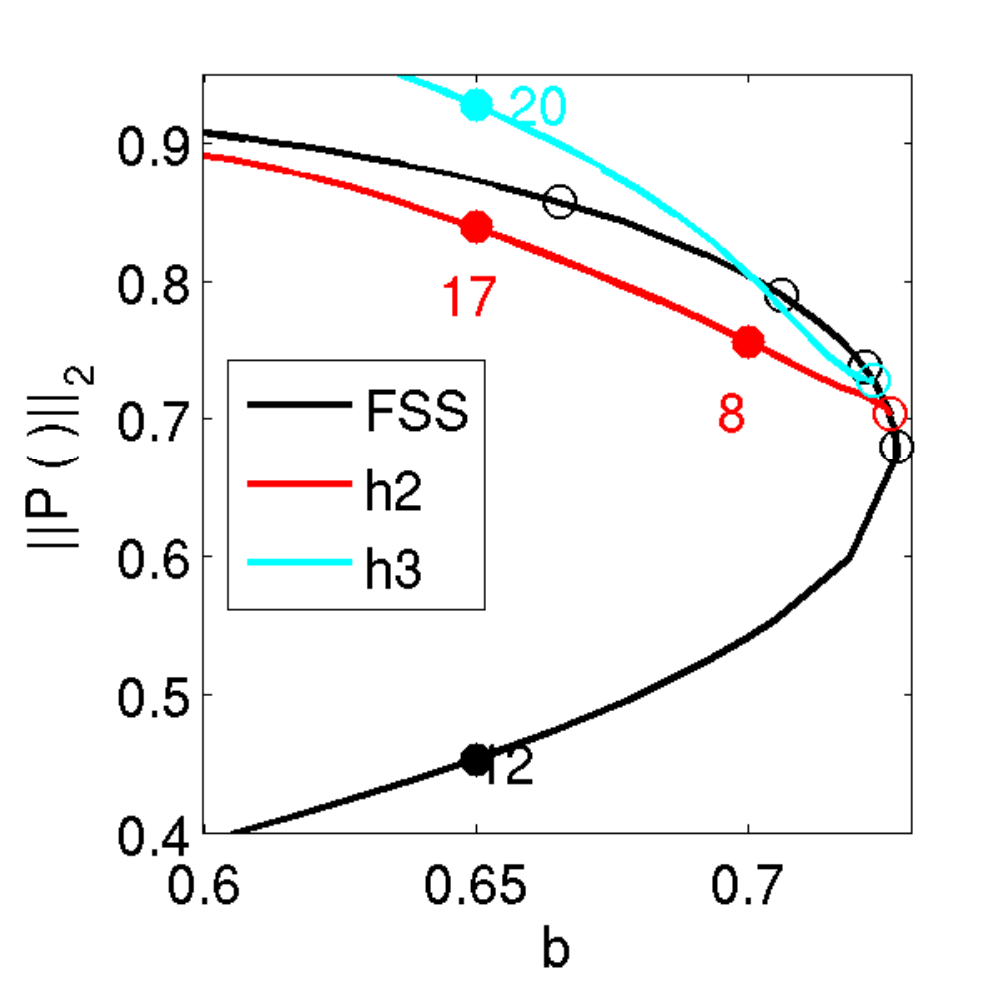}&
\raisebox{22mm}{\begin{tabular}{l}
\ig[width=38mm]{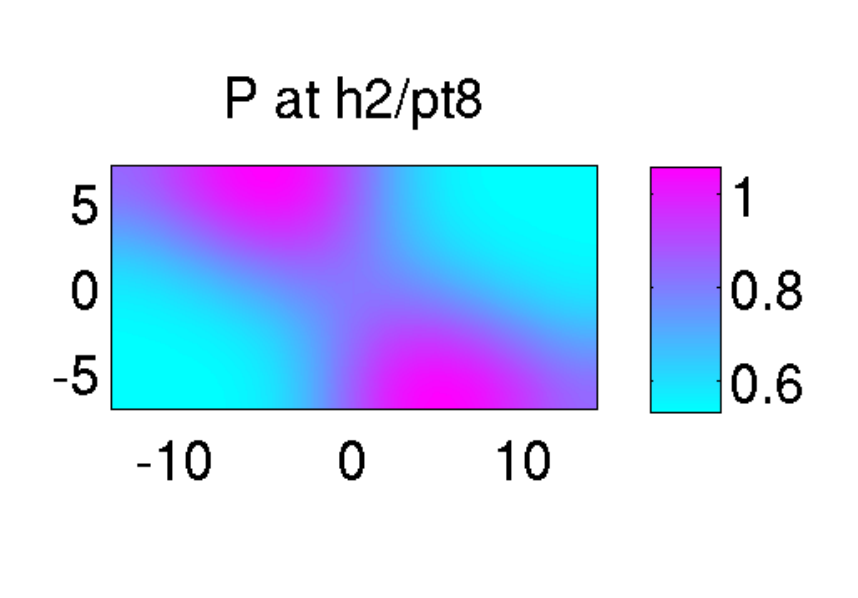}\ig[width=38mm]{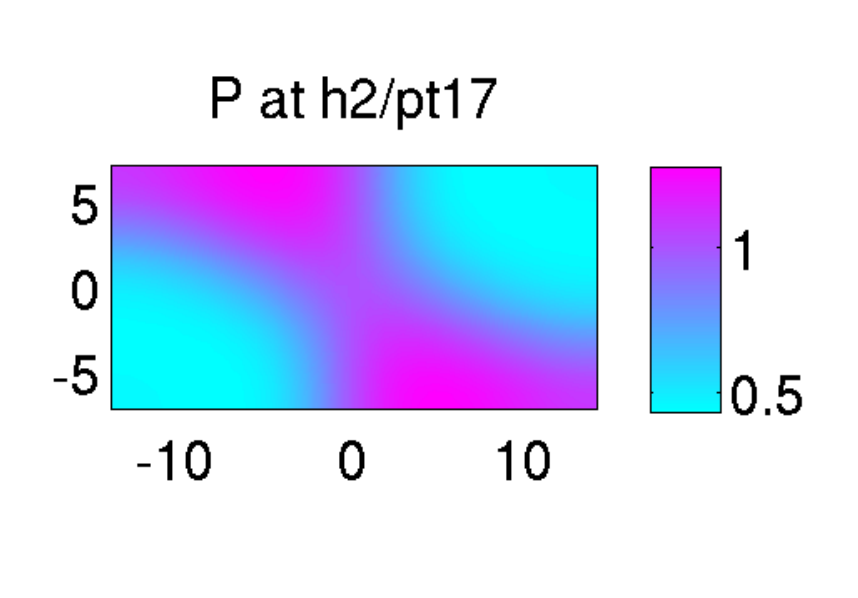}
\ig[width=38mm]{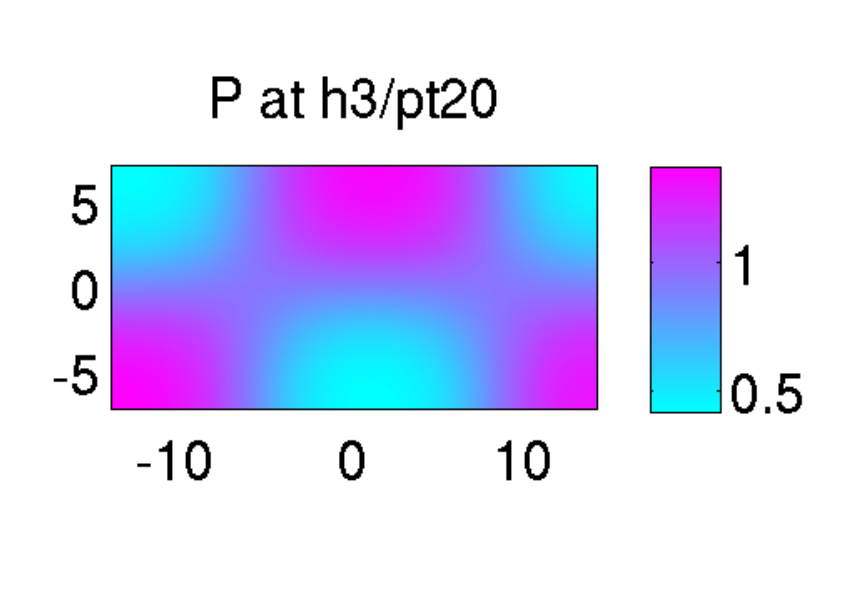}\\[-6mm]
\ig[width=38mm]{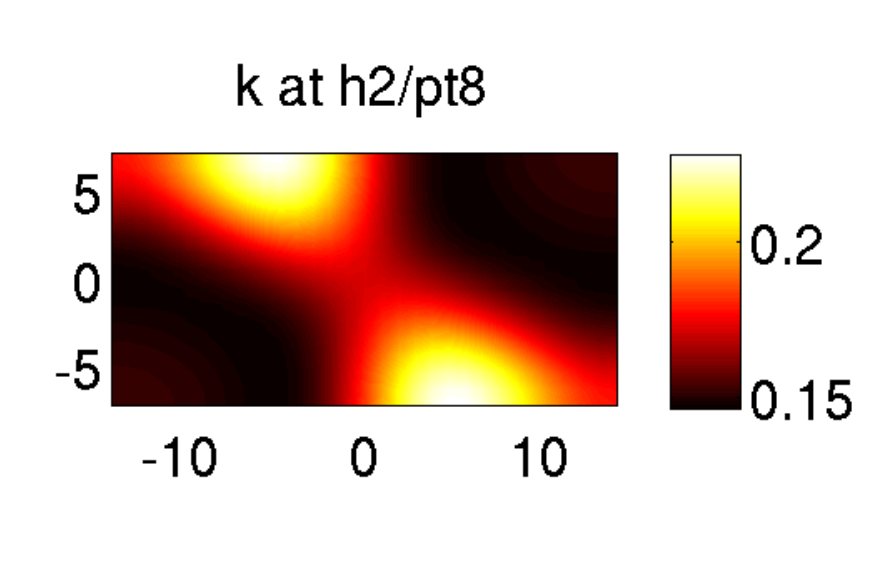}\raisebox{-0mm}{\ig[width=38mm]{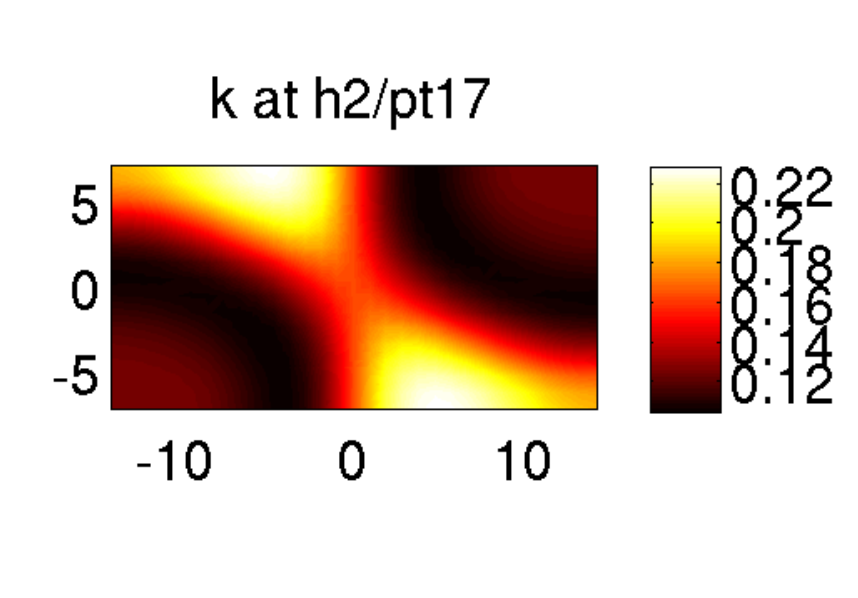}}
\ig[width=38mm]{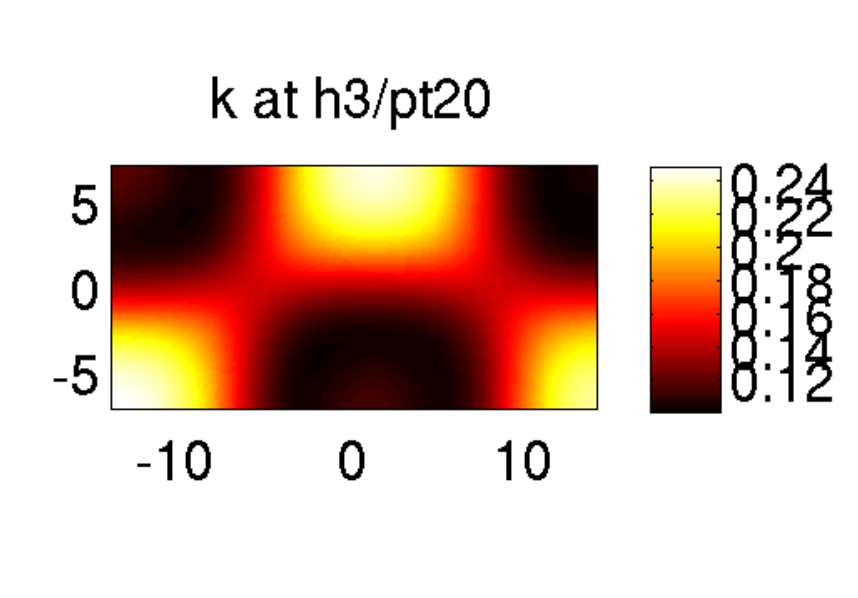}
\end{tabular}}
\end{tabular}
}
\caption{New patterned branches in 2D; $(x,y)\in(-L,L)\times(-\frac L 2, \frac L 2)$.  \label{2df1}}
\end{figure}

As a 2D example we consider \reff{slcan2} on the domain
$\Om=(-L,L)\times(-L/2,L/2)$, $L=2\pi/0.44$ as before, 
with a rather coarse mesh of $40\times
20$ points, hence approximately 1600 DoF.  The \fcss{} branches are of 
course the same
as in 1D (or 0D), and again at the end of the FSI branch we find a number
of Turing like bifurcations. In Fig.~\ref{2df1}(a),(b) we only present
the ``new'' patterned branches, i.e., those with a 
genuine $x$ and $y$ dependence. 

\begin{figure}[!ht]
\bce 
\ig[width=30mm,height=55mm]{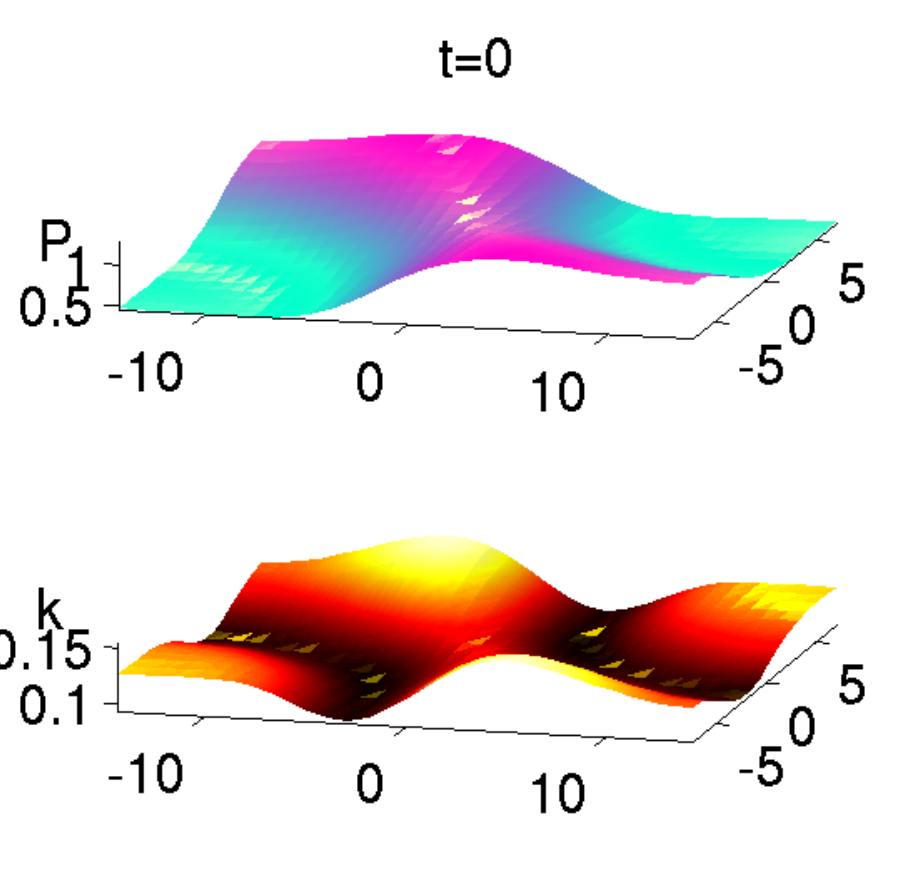}\quad\ig[width=30mm,height=55mm]{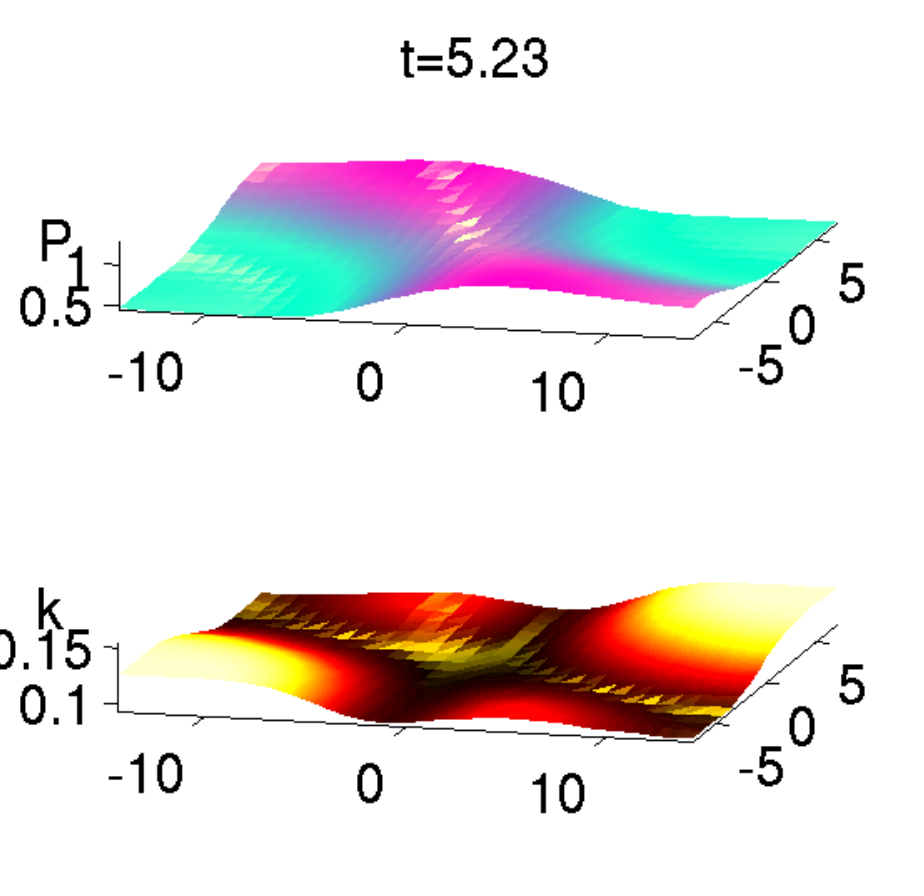}\quad
\ig[width=30mm,height=55mm]{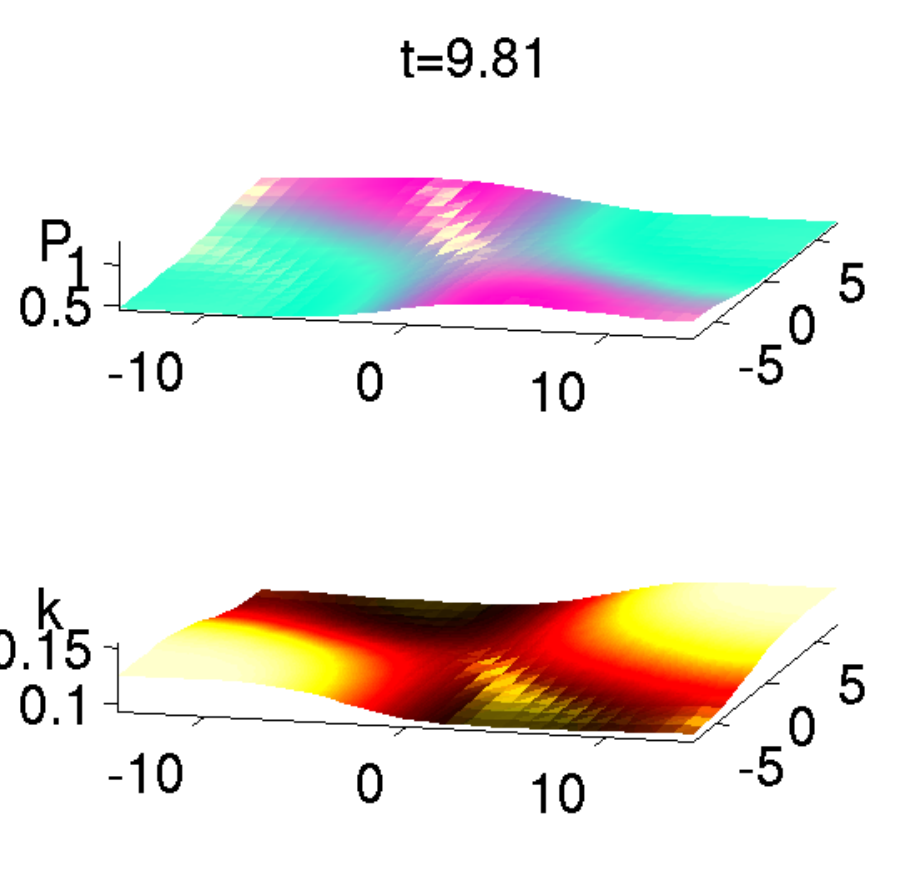}\quad\ig[width=30mm,height=55mm]{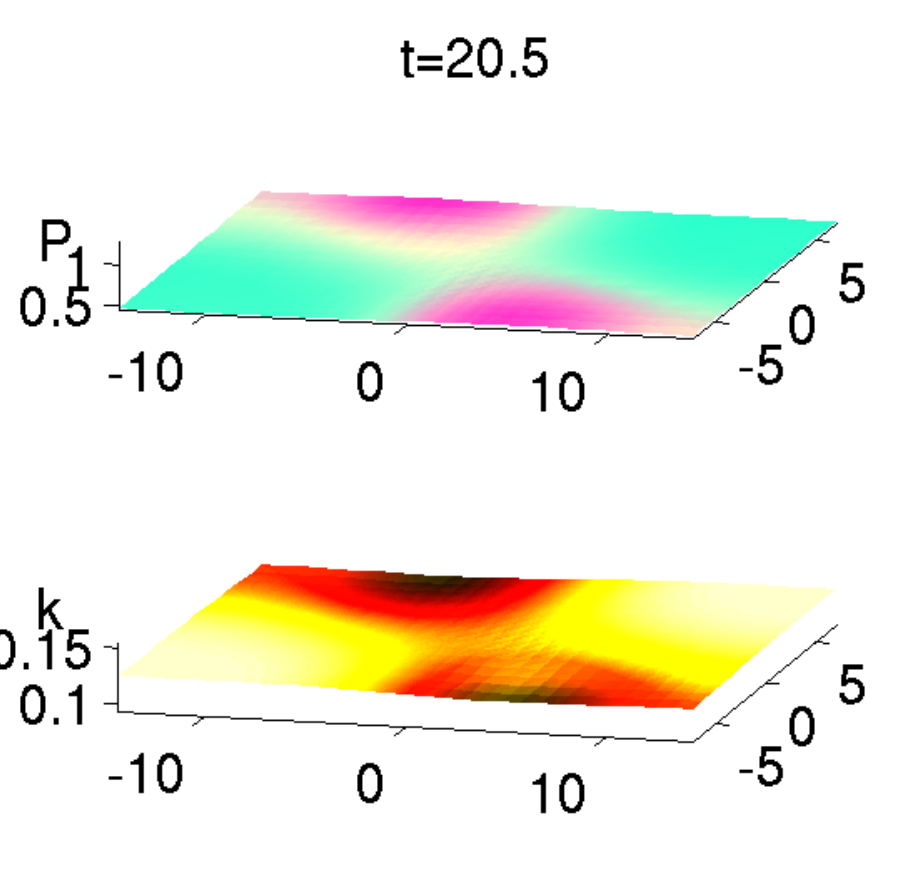}
\caption{Solutions on the canonical path to FSC\label{2df2}}
\ece 
\end{figure}

These new bifurcating \hcss{} again do not fulfill the SPP.  As an
example for a canonical path, in Fig.~\ref{2df2} we present 
snapshots from a path from $\hat{P}$ of 
the PCSS {\tt h2/pt17} to the FSC (see also 
\url{www.staff.uni-oldenburg.de/hannes.uecker/pde2path} for 
the movie), which yields a higher $J$ than the \hcss, i.e.,  
\begin{equation}\label{jord2}
  J({\rm \hcss})=-77.53<J({\rm \hcss}\to {{\rm FSC}})=-76.23<
J({\rm FSC})=-72.97. 
\end{equation}
Thus, this PCSS is not optimal, and neither is any other one we checked. 
Using the methods from \S\ref{1dsec} it is now of course also possible 
to find points with a genuine $x$ and $y$ dependence on the Skiba manifold 
separating FSC and FSM, but here we skip this presentation. 

The behavior and economic interpretation of the path from the PCSS 
to the FSC in Fig.~\ref{2df2} is 
rather similar to the convergence to the FSC in Fig.~\ref{ocf3}a. 
After a short transient the optimal strategy is to give a 
high phosphate load $k$ where $P$ 
is below the limit value $\hat P_{{\rm FSC}}$ (south-west and 
north-east corners of the domain), but initially there also is 
a high $k$ at high $P$ values (north-west and south-east corners). 

\section{Discussion}\label{dsec}
We have presented a numerical framework to treat infinite time horizon 
spatially distributed optimal control problems. First we derive 
the canonical PDE systems, which we then discretize in space and 
thus approximate by large systems of ODEs. For these we can resort to the theory 
and experience with non-distributed \OCPROs. Thus our results are 
intrinsically numerical; however, we believe that they can help to
develop the theoretical concepts for distributed \OCPROs, following 
Oliver Heavisides' saying 
{\em   ``Mathematics is an experimental science, and definitions do not
  come first, but later on.''}

From the economic point of view, the computation of canonical paths to the FOSS 
yields nontrivial and interesting results. Even more interesting would 
be locally stable POSS, but there is strong evidence that these 
do not exist for the shallow lake model \reff{sldiff1}, 
at least in the parameter regimes we considered so far \footnote{
see however ``Scenario 2'' in \cite{grass2014} for some locally stable POSS}.
On the other hand, in \cite{U15} we use our method to 
study the vegetation system from \cite{BX10}, and find 
that POSS dominate in large parameter regimes. We believe 
that the same happens in many other important systems, and a number 
of further investigations in this direction are under way.  
Natural candidates, i.e., systems with natural objective 
functions and controls, are 
related vegetation systems as in \cite{meron01, meron13}, 
fishery models as in \cite{Neu03, grass2012}, 
or ``crimo-taxis'' systems as in \cite{crimo10}.

\appendix
\section{SPP for PDEs}
\label{sec:spppdes}
In this appendix we discuss the SPP (Def.~\ref{def:spp}) in 
a somewhat more general situation, tailored to canonical systems 
coming from spatial discretizations of PDEs. 
Let $\uh=(\ph,\qh)\in \R^{2N}$ be a stationary state of a 
(non--distributed) canonical system of the form 
\hual{
\ddt \bpm p\\ q\epm&=F(p,q):=\bpm f(p,q)\\ 
rq-H_p(p,q))\epm, 
} 
where $f=H_q$, and let $J=D_u F(\uh)$ be the Jacobian at $\uh$. 
In \cite[Thm 7.10]{grassetal2008} it is explained 
that the eigenvalues of $J$ are symmetric around $r/2$, i.e., that 
there exist $N$ complex numbers $\xi_i$ such that 
\huga{\label{sJ}
\sigma(J)=\left\{\frac r 2 \pm \xi_i: i=1,\ldots, N\right\}. 
}
In detail, since $\det(J-\xi)=\det(J_r-(\xi-\frac r 2))$ where 
\huga{
J_r:=J-\frac r 2=\bpm H_{pq}-\frac r 2&H_{qq}\\
-H_{pp}&-H_{pq}+\frac r 2\epm, 
}
we have that $\frac r 2 +\xi_i\in\C$ is an eigenvalue of $J$ if and only if 
$\xi_i$ is an eigenvalue of $J_r$. But $J_r$ has the 
structure $\bpm A&B\\C&-A\epm $ with symmetric 
matrices $B,C\in\R^{N\times N}$, and as a consequence the eigenvalues of 
$J_r$ are $\xi_i=\pm \sqrt{\tilde\xi_i}$, $i=1,\ldots,N$. 

Now consider the distributed canonical system 
\hual{\label{pdecsys}
\pa_t \bpm p(x,t)\\ q(x,t)\epm &=F(p(x,t),q(x,t))
+\bpm D\Delta p(x,t)\\ -D\Delta q(x,t)\epm, 
} 
where $D\in\R^{N\times N}$ is a diffusion matrix, i.e., positive definite.  
Let 
\huga{\label{csdis} 
\ddt u(t)=G(u(t)), \quad u\in\R^{2nN}, 
}
be the associated spatially discretized system with $n$ spatial points, 
where 
$
u=(p(x_1),\ldots,p(x_n),\linebreak q(x_1),\ldots,q(x_n))\in\R^{2nN}, 
$
and let $\uh\in\R^{2nN}$ be a steady state of \reff{csdis}. 
Then $J=D_u G(\uh)$ has the structure 
$J=-K+\Jloc$, 
where $\Jloc$ has the block structure 
\huga{
\Jloc=\bpm H_{pq}^1&0&\ldots&0&H_{qq}^1&0&\ldots&0\\
0&H_{pq}^2&\ldots&0&0&H_{qq}^2&\ldots&0\\
\vdots&\vdots&\ddots&\ldots&\vdots&\vdots&\ddots&\vdots\\
0&0&\ldots&H_{qq}^n&0&0&\ldots&H_{qq}^n\quad\\[3mm]
-H_{pp}^1&0&\ldots&0&r-H_{qp}^1&0&\ldots&0\\
0&-H_{pp}^2&\ldots&0&0&r-H_{qp}^2&\ldots&0\\
\vdots&\vdots&\ddots&\ldots&\vdots&\vdots&\ddots&\vdots\\
0&0&\ldots&-H_{pp}^n&0&0&\ldots&r-H_{qp}^n
\epm, 
}
composed of local matrices $H_{pq}^j:=H_{pq}(x_j):=H_{pq}(p(x_j),q(x_j)), 
H_{qq}^j,H_{pp}^j\in\R^{N\times N}$, 
and $K=\bpm L&0\\ 0&-L\epm $ with $L\in\R^{nN}$ 
coming from the discretization of $D\Delta$. The notation $K$ 
of course reflects the FEM background of the present paper, but 
the same structure $\bpm L&0\\0&-L\epm$ occurs for any discretization, 
in any space dimension, and for any $D$ not necessarily diagonal, 
i.e., containing cross diffusion.  

It follows that again $\frac r 2+\xi_i$ 
is an eigenvalue of $J$ if and only if $\xi$ is an eigenvalue of 
$J_r:=J-\frac r 2$, where $J_{r}$ has the structure 
$$
J_{r}=\bpm A&B\\ C&-A\epm, \text{ with symmetric } B,C\in\R^{nN}.
$$
Applying \cite[Lemma B.2, Lemma B.3]{grassetal2008} we obtain 
\begin{theorem}\label{pthm1} Let $\uh$ be a steady state of the spatially discretized 
distributed system \reff{csdis}, and let $J$ be the associated Jacobian. 
Then there exist $\xi_i\in\C$, $i=1,\ldots, nN$, such that 
\huga{\label{sJ2}
\sigma(J)=\left\{\frac r 2 \pm \xi_i: i=1,\ldots, nN\right\}. 
}
\end{theorem} 
As a consequence, $\dim E_s(\uh)\le Nn$, and 
the only candidates $\uh$ for right BC in \reff{rbc} 
are those with the SPP. 
As a corollary we find a property that, on the discretized level, 
is equivalent to the SPP. 
\begin{cor} 
\label{prop:equivspp}
Let $\hat u\in\R^{2nN}$ be an equilibrium of the spatially discretized distributed system \reff{csdis} and $r>0$. Then $\hat u$ satisfies the SPP 
iff every eigenvalue $\eigval$ of the according Jacobian $J(\hat u)$ satisfies
\begin{equation}
\label{sppalt}
	\|\Re\eigval-\frac{r}{2}\|>\frac{r}{2}.
\end{equation}
\end{cor}

\begin{remark}\label{clrem}{\rm 
Theorem \ref{pthm1} is formulated on the discretized level, and one might 
ask how it ultimately relates to the PDE. As a first step one can ask: 
Let a steady state $\uh\in\R^{2nN}$ of \reff{csdis} be an approximation 
of a PDE steady state $(\hat p,\hat q)\in X$  for \reff{pdecsys},  
with $X\subset\{(p,q):\Om\ra \R^{2N}\}$ 
some function space, e.g., $X=[H^1(\Om)]^{2N}$. 
If $\tilde{\uh}\in\R^{2\tilde{n}N}$ is an approximation 
of $(\hat p,\hat q)$ on a finer mesh $\tilde{n}>n$, or just a different mesh, 
do we have 
\huga{\label{dimndep}
nN-\dim E_s(\uh)=\tilde{n}N-\dim E_s(\tilde{\uh})\quad ?
}

We do not want to go into the details here, but 
if $E_c(\uh)=\emptyset$, i.e., $\sigma(J)\cap\{{\rm Re} \xi=0\}=\emptyset$, 
then \reff{dimndep} 
is true, for large enough $n,\tilde n$. Given some $\uh$, this can be easily 
tested numerically, and it is also clear from 
an analytical point of view. Refining $\uh$ to 
$\tilde{\uh}$ we essentially add high frequency modes to 
the FEM (or finite difference) mesh. These introduce the same 
number of additional eigenvalues at large positive and negative $\xi$ for the 
linearization $\tilde{J}$, because $J_F(p,q):X\ra X$ is relatively 
compact with respect to the Laplacian, i.e., w.r.t.~
$(p,q)\mapsto (D\Delta p,-D\Delta q)$. On the other hand, 
the small eigenvalues $\mu_i$, $|\mu_i|<R$ for some fixed $R$, 
are only slightly perturbed, i.e., $|\mu_i-\tilde{\mu}_i|\le C 
\|\uh_*-\tilde{\uh}\|$, where $\uh_*$ is suitably 
defined, for instance by interpolating $\uh$ to the mesh of $\tilde{\uh}$. 
But $\|\uh_*-\tilde{\uh}\|\ra 0$ as $n,\tilde{n}\ra \infty$, 
which yields \reff{dimndep}. 

To make this rigorous, we need to define appropriate 
function spaces and study the approximation properties of the 
spatial discretization. This is easy, as the stationary problem for 
\reff{pdecsys} can 
be written as an elliptic system, and hence $(p,q)$ is arbitrary smooth, 
but we omit the details here. 

In fact, \reff{sppalt} can also be formulated on the PDE 
level and might therefore replace the SPP from Def.~\ref{def:spp} 
for spatially distributed models. 
However, we also postpone an in depth analysis of this to future work. 
}
\erem 

\bibliographystyle{alpha}
\newcommand{\etalchar}[1]{$^{#1}$}

\end{document}